\documentclass[10pt,french]{article}

\usepackage[francais]{babel}
\usepackage[latin1]{inputenc}
\usepackage[all]{xy}
\usepackage{amssymb}
\usepackage{pstricks}
\usepackage{amsmath}
\usepackage{array}
\usepackage{indentfirst}
\usepackage{hyperref}
\usepackage{graphicx}

\newcommand{\e}{\epsilon}
\newcommand{\F}{{\mathbb F}_2}
\newcommand{\id}{\mathrm{id}}
\def\M{{\mathbf M}}
\def\Mc{{\mathcal M}}
\newcommand{\m}{\mathrm{M}}
\def\L{{\mathbf L}}
\def\Sb{{\mathbf S}}
\def\Sc{{\mathcal S}}
\def\St{{\rm St}}
\def\Hom{{\rm Hom}}
\newcommand\End{{\rm End}}
\newcommand{\Ext}{{\rm Ext}}
\newcommand{\R}{\mathbb{R}}
\newcommand{\gl}{\mathrm{GL}}
\newcommand{\im}{\mathrm{im}}

\newcommand{\map}{\mathrm{map}}

\newcommand{\U}{\mathcal U}
\newcommand{\Z}{\mathbb Z}
\newcommand\ra{\rightarrow}

\newcommand\Sq{\mathrm{Sq}}
\newcommand{\A}{\mathcal A}
\newcommand{\diag}{\mathrm{diag}}
\newcommand{\pis}{\pi}
\newcommand{\reg}{\widetilde{reg}}

\newcommand\Ht{\widetilde{H}}
\newcommand{\proofof}[1]{{\noindent \it D\'emonstration de  #1 \quad}}
\newcommand{\bqn}{\begin{eqnarray*}}
\newcommand{\eqn}{\end{eqnarray*}}
\newcommand\qed{\hfill $\square$ \vspace{2ex}\par}

\author{Nguyen Dang Ho Hai, Lionel Schwartz, Tran Ngoc Nam
\footnote{The authors are partially supported by projet ANR blanc BLAN08-2\_338236, HGRT.}}
\title{La fonction g\'en\'eratrice de Minc et une  ``conjecture de Segal'' pour  certains spectres de Thom
}

\newtheorem{theorem}{Th\'eor\`eme}[section]
\newtheorem{proposition}[theorem]{Proposition}
\newtheorem{lemma}[theorem]{Lemme}
\newtheorem{corollary}[theorem]{Corollaire}

\newenvironment{remark}%
{
	\begin{trivlist}
	\refstepcounter{theorem}%
    	\item[]{\textbf{Remarque\ \thetheorem\ }}}%
    	{\end{trivlist}
}

\newenvironment{proof}
{\noindent\textit{D\'emonstration\quad}}{\hfill $\square$ \vspace{2ex}\par}

\begin{document}
\maketitle

\begin{abstract}
On construit dans cet article une r\'esolution injective minimale dans la
cat\'egorie $\U$ des modules instables sur l'alg\`ebre de Steenrod modulo $2$, de
la cohomologie de certains spectres obtenus  \`a partir de l'espace
de Thom du fibr\'e, associ\'e \`a la repr\'esentation r\'eguli\`ere r\'eduite du groupe
ab\'elien \'el\'ementaire $(\Z/2)^n$,  au dessus de l'espace $B(\Z/2)^n$.
Les termes de la r\'esolution sont des produits tensoriels de
modules de Brown-Gitler $J(k)$  et  de modules de
Steinberg $L_n$ introduits par S. Mitchell et S. Priddy. 
Ces modules sont injectifs d'apr\`es J. Lannes et S. Zarati, de plus ils
sont ind\'ecomposables. L'existence de cette r\'esolution avait \'et\'e
conjectur\'ee par Jean Lannes et le deuxi\`eme auteur. La principale
indication soutenant cette conjecture \'etait un r\'esultat
combinatoire de G. Andrews : la somme altern\'ee des s\'eries de
Poincar\'e des modules consid\'er\'ees est nulle.

Ce r\'esultat a des cons\'equences homotopiques et permet   
de d\'emontrer pour ces spectres un r\'esultat du type de la
conjecture de Segal pour les classifiants des 
$2$-groupes ab\'eliens \'el\'ementaires \cite{AGM85}.
\end{abstract}

\tableofcontents

%************************************************************************

\section{Introduction}\label{introduction}

Dans une cat\'egorie ab\'elienne il est en g\'en\'eral difficile de construire explicitement des r\'esolutions
injectives ou projectives minimales. C'est en particulier le cas dans la cat\'egorie des modules instables sur 
l'alg\`ebre de Steenrod modulo $2$ ${\mathcal A}$. 
On sait tr\`es bien  d\'ecrire les objets injectifs de la cat\'egorie \cite{LS89}, 
de plus comme ces modules  sont cohomologie modulo $2$ de spectres ou d'espaces 
(contrairement \`a ce qu'il en est pour les objets projectifs) ceci accroit 
l'int\'erêt pour de telles constructions. 
Cependant en dehors de quelques exemples et d'un r\'esultat de W. H. Lin \cite{Lin92}, 
peu utilisable, tr\`es peu de r\'esultats g\'en\'eraux sont connus. 
On n'a même pas de r\'esultats de finitude appropri\'e g\'en\'eral : par exemple si on sait que les modules 
ayant un nombre fini de g\'en\'erateurs ont des r\'esolutions dont chaque terme est somme directe finie 
d'injectifs ind\'ecomposables, on ne sait  pas  d\'emontrer l'analogue pour des modules instables dont 
l'enveloppe injective est elle même somme directe finie d'injectifs ind\'ecomposables 
(ce qui est la condition de finitude raisonnable pour la cohomologie modulo $2$ d'un espace de dimension infinie). 
Ce r\'esultat est \'equivalent \`a des conjectures difficiles concernant des cat\'egories de foncteurs entre espaces vectoriels 
sur le corps $\F$ (\cite{Dja07}).

Dans cet article on se propose d'\'etudier un cas sugg\'er\'e par certaines identit\'es combinatoires,  en fait on part
d'une formule montrant qu'une somme altern\'ee de s\'eries formelles est nulle. Dans la mesure o\`u \`a l'exception d'un 
terme les s\'eries formelles qui apparaissent sont les s\'eries de Poincar\'e de modules instables injectifs bien connus, 
que le terme restant est la s\'erie de Poincar\'e de la cohomologie d'un ``spectre de Thom'' on esp\`ere  r\'ealiser cette 
identit\'e alg\'ebriquement, c'est ce que nous faisons dans cet article,  puis g\'eom\'etriquement, ceci sera fait ailleurs. 
Mais le r\'esultat alg\'ebrique seul permet de d\'eduire des cons\'equences homotopiques, cela sera expliqu\'e plus bas.

La fonction de partition de Minc  $\nu(n)$ est d\'efinie comme le nombre de repr\'esentations de l'entier $n$ 
en somme d'entiers $c_i$ :  $n=c_1+\cdots +c_m$  avec $c_m \leq 2c_{m-1}\le \cdots \le 2^{m-1}c_1=2^{m-1}$,  $m$ quelconque.
On note $\nu(m,n)$ le nombre des solutions pour lesquelles $c_m \not =0$, $c_{m+1}=0$, $m$ donn\'e. On pose $\mu_m(q)= \sum_n \, \, \nu(m,n)q^n$. 
Dans \cite{And81} G. Andrews  montre que :
\begin{equation}\tag{\bf A}\label{andrew}
q^{2^n-1}\ell_m(q) = \sum_{i=0}^m(-1)^i \mu_{i}(q) \ell_{m-i}(q)
\end{equation}
avec
$$
\ell_m(q)=\frac{q^{(2-1)+(2^2-1)+\ldots+(2^m-1)}}{(1-q^{2-1})\ldots(1-q^{2^m-1})}.
$$

Soit $\U$ la cat\'egorie des modules instables sur l'alg\`ebre de Steenrod modulo $2$.
Les s\'eries formelles qui apparaissent ci-dessus sont celles du produit tensoriel de modules  
de Steinberg  $L_{m-i}$ \cite{MP83} et de Brown-Gitler $J(2^i-1)$ \cite{GLM92} qui sont 
des objets injectifs dans $\U$. Le module de Steinberg $L_n$ est un facteur direct dans 
$\F[x_1,\ldots,x_n]$, le module de Brown-Gitler $J(k)$ est lui un module fini caract\'eris\'e par
$\Hom_\U(M,J(k)) \cong M^{k*}$.
La s\'erie de Poincar\'e de $L_j$ est $\ell_j$, celle de $J(2^h-1)$ est $\mu_h$.

Le terme de gauche de l'\'egalit\'e (\ref{andrew}) 
est la s\'erie de Poincar\'e du sous-module   $L'_n = \omega_nL_n \subset L_n$. Ici $\omega_n$ est le produit de 
toutes les formes lin\'eaires non-nulles, c'est-\`a-dire la classe d'Euler de la somme de Whitney de  tous les fibr\'es 
en droites non triviaux sur $B(\Z/2)^n$. C'est la cohomologie d'un spectre de Thom  appropri\'e \cite{Tak99}
(voir \ref{cofibration-taka} ).
La s\'erie de Poincar\'e, $\ell'_n$, de $L'_n$
v\'erifie alors  $\ell'_n=t^{2^n-1}\ell_n$.
Pour un espace vectoriel gradu\'e $V$ on notera $P(V)$ sa s\'erie de Poincar\'e.
Le r\'esultat d'Andrews  dit que :

$$
-P(L'_n)+P(L_n)+\sum_{s=1}^n(-1)^s
P(L_{n-s}\otimes J(2^s-1))+(-1)^nP(J(2^n-1))=0.
$$
Ceci sugg\`ere la construction d'une r\'esolution injective pour  $L'_n$.
Voici le premier r\'esultat:

\begin{theorem}\label{NST1}
Pour tout $n\geq 1$, il
existe une r\'esolution injective minimale dans $\U$ :
$$ 0 \rightarrow  L'_n\rightarrow
L_n\to L_{n-1}\otimes J(1)\to L_{n-2}\otimes J(3)\rightarrow \cdots\rightarrow
L_1\otimes J(2^{n-1}-1)\rightarrow J(2^n-1)\rightarrow 0.
$$
\end{theorem}

On notera $f_{s,n}$ pour les
morphismes interm\'ediaires
$$
 L_{n-s+1}\otimes J(2^{s-1}-1)\to L_{n-s}\otimes J(2^s-1).
$$
On a  en corollaire :

\begin{theorem}
\label{NST2}
Soit $n\ge 2$ et soit le sous-ensemble de $\Z \times \Z$ d\'etermin\'e par :
$$\mathfrak A_n=\{t-s \le -2^{n-2}-n\}$$
On a $${\rm Ext}_{\mathcal A}^{s,t}(\Z/2,L'_n))=
\begin{cases}\Z/2 & \text{si $(s,t)\in \{ (n,1-2^n),(n+1,1-2^{n-1})\}$,}\\
0 & \text{si $(s,t)\in \mathfrak A_n\setminus \{(n,1-2^n),(n+1,1-2^{n-1})\}$.}
\end{cases}
$$
\end{theorem}

Ainsi qu'on l'a dit $L'_n$ est cohomologie modulo $2$ d'un spectre qui est obtenu comme suit. L'idempotent de Steinberg 
$e_n= \bar B_n \bar\Sigma_n \in \F[\gl_n]$ induit
une application sur $\Sigma B(\Z/2)^n$, le t\'elescope de cette application est (\`a suspension pr\`es)  
le spectre ${\bf M}(n)$ de cohomogie $M_n$ (voir ci-dessous). On peut aussi appliquer l'idempotent \`a 
l'espace de Thom du fibr\'e $\reg_n$ de base
$B(\Z/2)^n$ \cite{Tak99} qui est somme de tous les fibr\'es en droites non triviaux sur $B(\Z/2)^n$. 
On obtient alors comme t\'elescope de cette application  
(\`a suspension pr\`es c'est un espace) un spectre ${\bf L}(n)$ de cohomologie $L_n$. On peut encore appliquer cette proc\'edure au fibr\'e 
$\reg_n^{\oplus 2}$. On obtient alors comme t\'elescope de cette application  (\`a suspension pr\`es c'est un espace) 
un spectre ${\bf L}'(n)$ de cohomologie $L'_n$. Cette construction a été introduite par 
Shin-ishiro Takayasu \cite{Tak99} et sera détaillée en \ref{cofibration-taka}

Le th\'eor\`eme suivant a lieu pour ce spectre (\`a $2$-compl\'etion près), c'est un analogue de la conjecture de Segal 
(forme faible) \cite{AGM85,LZ87}:
\begin{theorem}
\label{NST3} 
\begin{enumerate}
\item
Pour $n\ge 2$, on a
$$\pis^{k}(\L'(n))=
\begin{cases}\Z/2 &\text{si $k\in \{ n+2^n-1, n+2^{n-1} \}$,}\\
0 &\text{si $k\in [n+2^{n-2},+\infty)\setminus \{n+2^n-1,n+2^{n-1}\}$.}
\end{cases}
$$
\item Pour le spectre $\L'(1)$, on a $\pis^1(\L'(1))=0$, 
$\pis^2(\L'(1))=\Z_2$ (l'anneau des entiers $2$-adique) et $\pis^k(\L'(1))=0$ si $k> 2$.
\end{enumerate}
\end{theorem}

En fait il semble clair que ces calculs peuvent être poussés plus loin, 
mais cela implique le calcul du foncteur division par les modules $L_n$ pour $n\ge 2$ sur l'algèbre de Dickson, 
ici on a seulement utilisé le cas $n=1$. Ceci sera étudié ailleurs.

L'article est organis\'e comme suit.
Dans la section \ref{injectifs}, on rappelle des r\'esultats concernant le facteur de 
Steinberg et les modules de Brown-Gitler. Dans la section \ref{exactness} on d\'emontre le th\'eor\`eme 
\ref{NST1} modulo une pr\'esentation de certains  modules de Brown-Gitler, celle ci est le coeur
de l'argument et est donn\'ee en section \ref{pres}.
Dans la section \ref{application-homology} 
on  donne des applications pour les groupes d'extensions et 
on d\'emontre le th\'eor\`eme \ref{NST2}. 
A l'aide de la suite spectrale d'Adams, on d\'emontre le th\'eor\`eme \ref{NST3} 
dans la section \ref{application-homotopy}. 

Les auteurs remercient le PICS Formath Vietnam du CNRS de les avoir soutenus
en facilitant leurs rencontres. Ils sont aussi reconnaissants \`a N. Kuhn pour ses commentaires utiles et pour avoir attir\'e leur attention vers 
la tour de Goodwillie de l'identit\'e \'evalu\'ee en sph\`eres impaires. 
Ils remercient \'egalement G. Powell de ses remarques judicieuses.

\section{Modules instables  injectifs}\label{injectifs}

Dans cette section on rappele ce qu'il convient sur les modules instables injectifs.

\subsection{Les modules de Steinberg}

Soit $\gl_n:=\gl_n(\F)$ le groupe des matrices $n\times n$ inversibles \`a
coefficients dans le corps \`a deux \'el\'ements $\F$. Ce groupe op\`ere \`a
gauche sur l'alg\`ebre polynomiale gradu\'ee $\F[x_1,\dots,x_n]$
(chaque g\'en\'erateur $x_i$ \'etant de degr\'e $1$) par la formule
\[
(\sigma\cdot
f)(x_1,\dots,x_n):=f(\sum_{j=1}^n\sigma_{j,1}x_j,\dots,\sum_{j=1}^n\sigma_{j,n}x_j),
\]
o\`u$\sigma=(\sigma_{i,j})_{n\times n}\in\gl_n$ et
$f\in\F[x_1,\dots,x_n]$. Cette action s'\'etend \'evidemment au
semi-groupe de toutes les matrices $\m_n(\F)$. L'alg\`ebre polynomiale
 $\F[x_1,\dots,x_n]$ est isomorphe \`a la cohomologie modulo $2$, $H^*(B(\Z/2)^n; \F)$, de
 l'espace $B(\Z/2)^n$. Cette cohomologie est un module instable
 sur $\A$, l'alg\`ebre de Steenrod modulo $2$, et les actions ci-dessus sont $\A$-lin\'eaires.

Soit $S$ un sous-ensemble du  groupe $\gl_n$. 
On note $\bar S
\in \F [\gl_n]$ la somme de
tous les \'el\'ements de $S$.
On consid\`ere les cas  du sous-groupe de Borel $B_n$
des matrices triangulaires sup\'erieures et du
sous-groupe $\Sigma_n$ des permutations. 
L'idempotent de Steinberg, $e_n$, est d\'efini par
par la formule $$e_n= \bar B_n \bar\Sigma_n.$$

\begin{proposition}[\cite{Ste56}] On a $e_n^2=e_n$ et le module $\F[\gl_n]e_n$
est projectif et absolument irr\'eductible.
\end{proposition}

 S. Mitchell et S. Priddy d\'efinissent le module de Steinberg \cite{MP83}  en
th\'eorie des modules instables par :
$$
M_{n}:=e_n\cdot\F[x_1,\dots,x_n].
$$
Comme
${\mathcal A}$  op\`ere de mani\`ere naturelle \`a gauche sur $\F[x_1,\dots,x_n]$
l'espace vectoriel $M_{n}$ est un sous-${\mathcal A}$-module de
$\F[x_1,\dots,x_n]$. D'apr\`es le th\'eor\`eme de  Carlsson-Miller \cite{Mil84}  
$\F[x_1,\dots,x_n]$ est injectif dans la cat\'egorie $\U$, comme $M_n$ en est un 
facteur direct il est \'egalement injectif.

On notera que dans \cite{MP83} l'action est \`a 
droite et que nous travaillons avec une action \`a gauche.
La version de
$M_{n}$ que nous utilisons n'est donc pas invariante par le groupe sym\'etrique 
mais par le sous-groupe de Borel $B_n$. La
proposition 2.6 de \cite{MP83} montre que quand applique les 
deux idempotents, $\bar B_n \bar\Sigma_n$ et $\bar\Sigma_n \bar B_n$,
\`a un $\A-\gl_n$-module  (module instable ayant une $\gl_n$- action compatible ) on obtient des modules instables isomorphes, 
les isomorphismes \'etant donn\'es par $\bar B_n$ et $ \bar \Sigma_n $.

L'alg\`ebre de Dickson $D(n)$ est l'alg\`ebre des \'el\'ements  
invariants sous l'action du groupe $\gl_n$ 
dans $\F[x_1,\ldots,x_n]$, elle est polyn\^omiale en des g\'en\'erateurs 
de degr\'e $2^n-2^{n-1}$,\ldots , $2^n-2^{i}$, \ldots, $2-1$.
Soit $\omega_n$  l'invariant de Dickson sup\'erieur en degr\'e $2^n-1$ : c'est le produit de toutes les formes lin\'eaires non-nulles, soit 
$$\omega_n=\det (x_j^{2^{i-1}})_{1\le i,j\le n},$$  
c'est aussi la classe d'Euler de la somme de tous les fibr\'es en droites r\'eelles non triviaux sur $B(\Z/2)^n$).

\begin{proposition}[\cite{MP83,Kuh87}] \label{dec}
Le module instable $M_{n}$ est un module sur l'alg\`ebre de Dickson $D(n)$. Le sous-espace vectoriel gradu\'e $L_n=\omega_nM_n \subset M_n$ est un sous module instable.
De plus il y a un isomorphisme de ${\mathcal A}$-modules  : $$ M_{n} \cong L_n
\oplus L_{n-1}.$$ Cet isomorphisme est
rigide.
\end{proposition}

La d\'ecomposition des modules instables $M_{n-1} \cong L_n
\oplus L_{n-1}$ correspond \`a une d\'ecomposition $e_n =
\e_n+\e'_n$ dans l'alg\`ebre du semi-groupe des matrices
$\F[M_n(\F)]$. Les \'el\'ements $\e_n$ et $\e'_n$ sont des idempotents
primitifs et orthogonaux, $\e'_n$ est constitu\'e de matrices
singuli\`eres.

Par d\'efinition, $M_1$ est la cohomologie $H^*(B\Z/2;\F)$. On en d\'eduit que 
$L_1$ est la cohomologie r\'eduite $\Ht(B\Z/2;\F)$. 
En identifiant $L_1^{\otimes n}$ \`a l'id\'eal de $\F[x_1,\cdots,x_n]$ engendr\'e par $x_1\cdots x_n$, on 
peut v\'erifier que

\begin{proposition} \label{cap-steinberg}
$L_n=e_n\cdot L_1^{\otimes n}=\omega_nM_{n}=\bigcap_{i=1}^{n-1} L_1^{\otimes i-1}\otimes L_2 \otimes L_1^{\otimes n-i-1}.$
\end{proposition}

Celle-ci sera d\'emontr\'ee en appendice en utilisant la relation qui existe entre les idempotents de Steinberg et 
l'alg\`ebre de Hecke $\End_{\F[\gl_n]}(1_{B_n}^{\gl_n})$ \'etudi\'ee par N. Kuhn \cite{Kuh84}.

Pour tout $1\le k\le n$, l'inclusion  
$$\bigcap_{i=1}^{n-1} L_1^{\otimes i-1}\otimes L_2 \otimes L_1^{\otimes n-i-1}\subset 
\big(\bigcap_{i=1}^{k-1} L_1^{\otimes i-1}\otimes L_2 \otimes L_1^{\otimes n-i-1}\big) \cap \big( 
\bigcap_{i=k+1}^{n-1} L_1^{\otimes i-1}\otimes L_2 \otimes L_1^{\otimes n-i-1}\big )$$
d\'efinit une inclusion canonique $\delta\colon L_n\hookrightarrow L_k\otimes L_{n-k}$. Il 
est clair que $\delta$ est coassociative. On obtient donc une structure de $\F$-coalg\`ebre 
sur $L_*:=\bigoplus_{i\ge 0}L_i$ qui m\'eritera une \'etude ailleurs.

La s\'erie de Poincar\'e d'un espace vectoriel gradu\'e $V$ est d\'efinie
par $P(V)=P(V,t):=\sum_d\dim V^dt^d$, o\`u $V^d$ d\'esigne la partie
de degr\'e $d$ de $V$. On consid\`ere aussi la s\'erie de Poincar\'e de l'espace
vectoriel sous-jacent d'un module instable $M$, on la notera aussi
$P(M)$ par abus. On a :

\begin{proposition}[\cite{MP83}] 
La s\'erie de Poincar\'e de $L_n$, not\'ee $\ell_n$, est
donn\'ee par
$$
\ell_n=\prod_{i=1}^n\frac{t^{2^i-1}}{1-t^{2^i-1}}.
$$
\end{proposition}

En fait Mitchell et Priddy montrent qu'en tant qu'espace vectoriel gradu\'e,
$\bar\Sigma_n \bar B_n \cdot \F[x_1,\dots,x_n]$ a une base
form\'ee par les \'el\'ements  $$\Sq^{i_1+1}\cdots \Sq^{i_n+1}(\frac{1}{x_1\cdots x_n})$$
o\`u$i_1> 2i_2>\cdots> 2^{n-1}i_n\geq 0$. 
La copie de $M_n$ que l'on consid\`ere a donc pour base les \'el\'ements
$$\bar B_n \cdot \Sq^{i_1+1}\cdots \Sq^{i_n+1}(\frac{1}{x_1\cdots x_n}).$$

\begin{theorem}\label{base-steinberg}
En tant que $\F$-espace vectoriel gradu\'e,
le module $M_{n}$ a  une base form\'ee par les \'el\'ements
$$e_n \cdot \omega_1^{i_1-2i_2}\cdots \omega_{n-1}^{i_{n-1}-2i_n} \omega_n^{i_n},$$
o\`u $i_1> 2i_2>\cdots> 2^{n-1}i_n\geq 0.$ Ici $\omega_k\in  \F[x_1,\cdots,x_k]$.
\end{theorem}
Ceci sera d\'emontr\'e en appendice.
On notera que l\'el\'ement  $e_n \cdot \omega_1^{i_1-2i_2}\cdots \omega_{n-1}^{i_{n-1}-2i_n} \omega_n^{i_n}$ 
est de degr\'e $i_1+\cdots+i_n$.

\subsection{Les modules de Brown-Gitler}

Soit $J(k)$ le ${\mathcal A}$-module de Brown-Gitler (cf. \cite[Chapter 2]{Sch94}).  
En degr\'e $k$, l'espace vectoriel gradu\'e $J(k)$ est \'egal \`a $\F$, engendr\'e par un
\'el\'ement not\'e $\iota_k$. Le module $J(k)$ est caract\'eris\'e par le fait que 
la transformation naturelle qui \`a $f \in {\mathrm{Hom}}_\U(M, J(k))$ associe sa restriction
en degr\'e $k$, qui est donc dans le dual $M^{k*}$, est une \'equivalence naturelle :
$$
{\rm Hom}_\U(M,J(k) \cong M^{k*}.
$$
En particulier
 si un ${\mathcal A}$-module instable $M$ est de dimension $1$ en degr\'e $k$,
alors il existe un et un seul morphisme ${\mathcal A}$-lin\'eaire non nul de
degr\'e z\'ero de $M$ dans $J(k)$; ce morphisme envoie sur $\iota_k$
l'\'el\'ement non nul de $M^k$ (la partie de degr\'e $k$ de $M$).

H. Miller a donn\'e dans \cite{Mil84} une description globale des $J(k)$ en consid\'erant leur somme directe. 
Il introduit l'objet bigradu\'e $J^*_*$ d\'etermin\'e par $J^\ell_k=J(k)^\ell$. 
Cet objet est en fait une alg\`ebre bigradu\'ee, dot\'ee d'une structure de module instable pour laquelle la formule de Cartan a lieu.
En fait Miller d\'emontre que (\cite[Chapter 2]{Sch94}) :

\begin{proposition}\label{J**} On a
$$
J^*_* \cong \F[\hat t_i\mid i\geq 0]
$$
avec $\hat t_i \in J(2^i)^1$ de bidegr\'e $(1,2^i)$. La structure de
${\mathcal A}$-module instable de $J_*^*$ est d\'etermin\'ee par
$$\Sq^1(\hat t_i)=\hat t_{i-1}^2, \quad i\geq 1, \quad \Sq^1 (\hat t_0)=0$$
et la formule de Cartan. 
Le module $J(k)$ s'identifie au sous-espace engendr\'e par les mon\^omes de second degr\'e $k$, {\it i.e.}  
par les $\hat t_0^{\alpha_0}\cdots \hat t_i^{\alpha_i}$
avec $\sum_h \alpha_h2^h=k$.
\end{proposition}

Soit maintenant $\Omega_k$ l'ensemble des suites d'entiers $(i_1,\dots,i_k)$
telles que
$$
0<i_1\le 2i_2\le 4i_3\le \cdots \le 2^{k-1}i_k=2^{k-1}.
$$
 La $k$-i\`eme fonction g\'en\'eratrice de Minc   $\mu_k$, est  donn\'ee par
$$\mu_k(t)=\sum_{d\geq 0}
{|\Omega_k^d|}t^d=\sum_{\Omega_k}\; t^{i_1+\cdots+i_k},$$
$\Omega_k^d$ \'etant le sous ensemble de $\Omega_k$ constitu\'e par les partitions de somme $d$.

\begin{proposition}[{\cite[p. 57]{Sch94}}] 
\label{poincare-J}
Soit $k \geq 1$, on  a $P(J(2^k-1))=\mu_k.$
\end{proposition}

Dans la r\'ef\'erence ceci est propos\'e  en exercice. La d\'emonstration r\'esulte de \ref{J**}.
Partant du mon\^ome $\hat t_0^{\alpha_0}\cdots \hat
t_{k-1}^{\alpha_{k-1}}\in J(2^k-1)^d$,
on pose $$i_1=\frac{\alpha_0}{2}+\frac{1}{2}, i_{2}=\frac{\alpha_1}{2}+\frac{\alpha_0}{4}+\frac{1}{4},\ldots, i_k=\frac{\alpha_{k-1}}{2}+\cdots +\frac{\alpha_0}{2^k}+\frac{1}{2^k}=1.$$ On v\'erifie facilement
que les $i_h$ sont des entiers, et que la suite $(i_1,\ldots,i_k)$ est dans $\Omega_k^d$. Inversement partant d'une suite $(i_1,\ldots,i_k) \in \Omega_k^d$, les formules
$$\alpha_0=2i_1-1, \alpha_1=2i_2-i_1,\dots,\alpha_{k-1}=2i_k-i_{k-1}$$ 
d\'eterminent un mon\^ome comme ci-dessus, fournissant l'application r\'eciproque. Le r\'esultat suit.

\subsection{Le th\'eor\`eme de Lannes-Zarati}

Enfin on rappelle que:

\begin{theorem}
Le module instable $L_h \otimes J(k)$ est injectif dans la cat\'egorie $\U$.
\end{theorem}

C'est un cas particulier du r\'esultat principal de Lannes et Zarati dans \cite{LZ86}. Par ailleurs,
il r\'esulte de \cite{LS89} que ce module est ind\'ecomposable.

%****************************************************************

\section{Construction des morphismes et exactitude}\label{exactness}

Dans cette section on construit les morphismes du complexe, puis on  d\'emontre l'exactitude,  
modulo une pr\'esentation du module de Brown-Gitler $J(2^k-1)$ qui sera faite dans la section suivante.

\subsection{Construction des morphismes}

Rappelons que $L_1$, \'etant la cohomologie r\'eduite de $B\Z/2$, s'identifie \`a l'id\'eal $(x)\subset \F[x]$. 
On note $\pi_s\colon L_1\rightarrow J(2^s)$ l'unique morphisme non trivial qui envoie $x^{2^s}$ 
sur la classe fondamentale $\iota_{2^s}$ de $J(2^s)$.

D\'efinissons le morphisme $f_{s,n}$ comme suit:
$$\xymatrix
{
L_{n-s+1}\otimes J(2^{s-1}-1)\ar[rr]^-{f_{s,n}}\ar[d]^-{\delta\otimes \id}	&	&	 L_{n-s}\otimes J(2^s-1)\\
L_{n-s}\otimes L_1\otimes J(2^{s-1}-1)\ar[rr]^-{\id\otimes \pi_{s-1}\otimes \id} &  &
L_{n-s}\otimes J(2^{s-1})\otimes J(2^{s-1}-1)\ar[u]^-{\id\otimes \mu}
}
$$
Ici $\mu\colon J(2^{s-1})\otimes J(2^{s-1}-1)\rightarrow J(2^s-1)$ est 
la multiplication, qui est l'unique morphisme non trivial et $\delta\colon L_{n-s+1}\rightarrow L_{n-s}\otimes L_1$ 
la comultiplication de la $\F$-coalg\`ebre $L_*$.
Par convention, l'inclusion naturelle $L'_n\hookrightarrow L_n$ se note $f_{0,n}$.

\begin{proposition} $f_{s+1,n}\circ f_{s,n}=0$ pour  $1\le s\le n-1$.
\end{proposition}
\begin{proof}On va se ramener au cas $n=2$. 
Grâce \`a la coassociativit\'e de $\delta$ et \`a l'associativit\'e de $\mu$,  
la compos\'ee $f_{s+1,n}\circ f_{s,n}$, pour $1\le s\le n-1$, 
se factorise alors comme suit:
$$\xymatrix
{
L_{n-s+1}\otimes J(2^{s-1}-1)\ar[r]^{\delta\otimes \id}\ar[d]^{\delta\otimes \id}\ar@/_9pc/[dd]|-{f_{s+1,n}\circ f_{s,n}}
	&	L_{n-s-1}\otimes L_2\otimes J(2^{s-1}-1)\ar[d]|{\id \otimes \delta\otimes \id}\\
L_{n-s}\otimes L_1\otimes J(2^{s-1}-1)\ar[r]^-{\delta\otimes \id\otimes \id}	&	L_{n-s-1}\otimes L_1\otimes L_1 \otimes J(2^{s-1}-1)\ar[d]|-{\id \otimes \pi_s\otimes \pi_{s-1}\otimes \id}\\
L_{n-s-1}\otimes J(2^{s+1}-1)		&		L_{n-s-1}\otimes J(2^{s})\otimes J(2^{s-1}) \otimes J(2^{s-1}-1).\ar[l]_-{\id\otimes \mu}
}
$$
A cause du corollaire \ref{trivial} ci-dessous, on a   $ \mu \circ  \circ \pi_s \otimes \pi_{s-1} \circ \delta =0$, ici $\mu$ d\'esigne la multiplication $J(2^{s})\otimes J(2^{s-1} \ra J(2^s +2^{s-1})$. Et donc $f_{s+1,n}\circ f_{s,n}=0$ pour  $1\le s\le n-1$. 
\end{proof}

On part de la base
comme espace vectoriel gradu\'e du facteur $L_2$
constitu\'ee  par les \'el\'ements
$e_2 \cdot \omega_1^{a-2b}\omega_2^b$ avec $a>2b>0$.

\begin{lemma}
Si $a>2b>0$ et $a+b=2^i+2^{i-1}$, alors l'expression de $e_2 \cdot \omega_1^{a-2b}\omega_2^b$
comme somme de mon\^omes distincts
ne contient pas $x_1^{2^i}x_2^{2^{i-1}}$.
\end{lemma}

\begin{proof} Notons que les conditions du lemme impliquent que $i>1$.
On a
\begin{equation*} \begin{split}
e_2 \cdot \omega_1^{a-2b}\omega_2^b
&= [x_2^{a-2b}+(x_1+x_2)^{a-2b}]x_1^bx_2^b(x_1+x_2)^b\\
&= \sum_{j=1}^b(\binom{b}{j}+\binom{a-b}{j})x_1^{b+j}x_2^{a-j}+
\sum_{j=b+1}^{a-b}\binom{a-b}{j}x_1^{b+j}x_2^{a-j}.
\end{split}\end{equation*}
Comme $a>2b$ et $a+b=2^i+2^{i-1}$, on voit que $a>2^i$ et $b<2^{i-1}$.
Posons $a=2^i+c$ et $b=2^{i-1}-c$ avec $0< c < 2^{i-1}$.
On en d\'eduit que le premier terme dans la somme ci-dessus ne peut contenir
$x_1^{2^i}x_2^{2^{i-1}}$.
D'autre part,  le coefficient de $x_1^{2^i}x_2^{2^{i-1}}$ dans le deuxi\`eme terme est
$\binom{a-b}{2^i-b}$.
Supposons que $c=2^tc^\prime$ avec $0\le t<i-1$ et $c^\prime$ impair.
On a alors
$$\binom{a-b}{2^i-b}=\binom{2^{i-1}+2c}{2^{i-1}+c}=\binom{2^{i-1}+2c}{c}=
\binom{2^{i-1}+2^{t+1}c^\prime}{2^tc^\prime}=\binom{2^{i-1-t}+2c}{c^\prime}=0$$
car $2^{i-1-t}+2c$ est pair alors que $c^\prime$ est impair.
\end{proof}

Il suit :

\begin{corollary}\label{trivial}
La compos\'ee
$$L_2\hookrightarrow L_1\otimes L_1 \xrightarrow{\pi_i\otimes \pi_{i-1}} J(2^i)
\otimes J(2^{i-1})\xrightarrow{\mu} J(2^i+2^{i-1})$$
est nulle pour tout $i\geq 1$.
\end{corollary}

\begin{proof}Le cas $i>1$ vient du lemme pr\'ec\'edent.
Si $i=1$, on v\'erifie que $L_2$ est trivial en degr\'es inf\'erieurs \`a $4$.
\end{proof}

\subsection{D\'emonstration de l'exactitude}
On commence par introduire l'application suivante :

$$g_s\colon L_1^{\otimes s} \xrightarrow{\pi_{s-1}\otimes \cdots\otimes
\pi_0}J(2^{s-1})\otimes \cdots \otimes J(1)\stackrel{\mu}\longrightarrow J(2^s-1)$$
o\`u$\mu$ est l'unique application non-triviale.
On montrera dans la section suivante que :

\begin{proposition}
L'application $g_s$ est surjective. Le syst\`eme d'\'el\'ements
$g_s(x_1^{i_1}\cdots x_s^{i_s})$ avec
$(i_1,\dots,i_s)\in\Omega_s$ est  une base de
$J(2^s-1)$.
\end{proposition}

Afin d'all\'eger les notations, on notera $\omega^{i_1,\ldots,i_{n-s}}$ l'\'el\'ement
$$e_{n-s}\cdot \omega_1^{i_1-2i_2}\cdots
\omega_{n-s-1}^{i_{n-s-1}-2i_{n-s}} \omega_{n-s}^{i_{n-s}}$$
et $g^{i_{n-s+1},\ldots,i_n}$ l'\'el\'ement $$g_s(x_{n-s+1}^{i_{n-s+1}}\cdots x_n^{i_n}).$$
La proposition suivante est la cons\'equence de la proposition pr\'ec\'edente et de \ref{base-steinberg}.

\begin{proposition} Soit  $0\le s\le n$. En tant qu'espace vectoriel gradu\'e,
$L_{n-s}\otimes J(2^s-1)$ a une base form\'ee par les \'el\'ements
\[\omega^{i_1,\ldots,i_{n-s}}\otimes g^{i_{n-s+1},\ldots,i_n},
\]
avec $i_1> 2i_2>\cdots> 2^{n-s-1}i_{n-s}> 0 <i_{n-s+1}\le 2i_{n-s+2}\le\cdots
 \le 2^{s-1}i_n=2^{s-1}.$
\end{proposition}

Pour $1\le s\le n$, notons $A(s,d)$ l'ensemble des \'el\'ements de cette base
qui v\'erifient $i_1+\cdots+i_n=d$ et $i_{n-s}\le 2i_{n-s+1}$, $B(s,d)$
l'ensemble de ceux qui v\'erifient $i_1+\cdots+i_n=d$ et
$i_{n-s}> 2i_{n-s+1}$.

Pour $s=0$, notons $A(0,d)$ l'ensemble des
$\omega^{i_1,\ldots,i_n}$ avec $i_1+\cdots+i_n=d$ et
$i_1> 2i_2>\cdots> 2^{n-1}i_{n}=2^{n-1}$, $B(0,d)$
l'ensemble des
$\omega^{i_1,\ldots,i_n}$ avec $i_1+\cdots+i_n=d$ et
$i_1> 2i_2>\cdots> 2^{n-1}i_{n}>2^{n-1}$.

Il est
clair qu'en degr\'e $d$, la dimension de $L_{n-s}\otimes J(2^s-1)$
est $|A(s,d)|+|B(s,d)|$, la somme des cardinaux de
$A(s,d)$ et $B(s,d)$. De plus $|A(s,d)|=|B(s+1,d)|$.

\begin{lemma}\label{diminfo}
Soit $1\le s\le n$. Alors en degr\'e $d$, $\dim\im f_{s,n}\geq |A(s-1,d)|$.
\end{lemma}
\begin{proof} D'apr\`es \ref{diagonal}, on a
$$\omega^{i_1,\ldots,i_{n-s+1}}=\omega^{i_1,\ldots,i_{n-s}}\cdot x_{n-s+1}^{i_{n-s+1}}+
\sum_{i}f_i \cdot x_{n-s+1}^i$$
pour certains
$i>i_{n-s+1}$ et $f_i \in \F[x_1,\ldots,x_{n-s}]\cdot\omega_{n-s}$.
D'où
$$f_{s,n}(\omega^{i_1,\ldots,i_{n-s+1}}\otimes g^{i_{n-s+2},\ldots,i_n})=
\omega^{i_1,\ldots,i_{n-s}}\otimes g^{i_{n-s+1},\ldots,i_n} +\sum_{i}y_i\otimes z_i
$$
pour certains $y_i\in L_{n-s}$ et $z_i\in J(2^s-1)$.
Comme $\deg z_i>\deg g^{i_{n-s+1},\ldots,i_n}$,
il suit facilement de cette formule que les
\'el\'ements de $f_{s,n}(A(s-1,d))$ sont lin\'eairement ind\'ependents.
\end{proof}

On d\'emontre l'exactitude de la suite dans le th\'eor\`eme \ref{NST1}:
$$ 0 \rightarrow  L'_n\rightarrow
L_n\xrightarrow{f_{1,n}} L_{n-1}\otimes J(1)\to L_{n-2}\otimes
J(3)\rightarrow \cdots\rightarrow
L_1\otimes J(2^{n-1}-1)\xrightarrow{f_{n,n}} J(2^n-1)\rightarrow 0.
$$
L'exactitude en $L_n$ est facile de v\'erifier en utilisant
\ref{diagonal} et la d\'efinition de $L'_n$.
L'exactitude en $J(2^n-1)$, {\it i.e.} la surjectivit\'e de $f_{n,n}$,
sera montr\'ee dans la section \ref{pres}.
Soit  $1\le s\le n-1$.
D'apr\`es \ref{diminfo}, en tout degr\'e $d$, on
a
\bqn \dim\im f_{s,n}+\dim\im
f_{s+1,n}&\geq&|B(s,d)|+|A(s,d)|\\ &=&\dim
L_{n-s}\otimes J(2^s-1). \eqn
Comme $\im f_{s,n}\subset\ker
f_{s+1,n}$, il suit que cette in\'egalit\'e est en fait une \'egalit\'e et donc ${\dim \im}f_{s,n}={\dim \ker}f_{s+1,n}$.
Cela prouve l'exactitude en $L_{n-s}\otimes J(2^s-1)$ pour
$1\le s\le n-1$. Le r\'esultat suit.

\section{Une pr\'esentation de $J(2^n-1)$}\label{pres}

\subsection{La pr\'esentation de $J(2^n-1)$}

Dans cette section on donne une description de $J(2^n-1)$ comme quotient de
l'id\'eal $(x_1\cdots x_n)\subset \F[x_1,\ldots,x_n]$.
D'apr\`es \ref{J**}   $J(2^n-1)$ est  le sous-module de $J_*^*$ qui admet pour base les mon\^omes de second degr\'e $2^n-1$, 
c'est-\`a-dire les mon\^omes $\hat t_0^{\alpha_0}\ldots \hat t_{k-1}^{\alpha_{k-1}}$ tels que $\sum_{h=0}^{n-1}\alpha_h2^h=2^n-1$.

On d\'esignera par
$MP(i)$ le sous-module $L_1^{\otimes i-1}\otimes L_2\otimes L_1^{\otimes n-i-1}$,
$1\le i\le n-1$ et par $MP(n)$ le sous-module $L_1^{\otimes n-1}\otimes L'_1$.
On consid\`ere donc :
$$g_n\colon L_1^{\otimes n}\xrightarrow{\pi_{n-1}\otimes \cdots\otimes
\pi_0}  J(2^{n-1})\otimes \cdots \otimes J(1)\stackrel{\mu}\longrightarrow J(2^n-1).$$
Par \ref{trivial} et le fait que $\pi_0(L'_1)$ est trivial,
le noyau de $g_n$ contient la somme $MP(1)+\cdots+MP(n)$.

\begin{theorem}\label{presentation} L'application $g_n$ est surjective et induit un isomorphisme
de modules instables
$$\frac{L_1^{\otimes n}}{MP(1)+\cdots+MP(n)} \cong J(2^n-1).$$
\end{theorem}

\begin{proof}En supposant que  $g_n$ est surjectif, la d\'emonstration de
l'isomorphisme se fait comme suit.
On rappelle que $\Omega_n$ est l'ensemble des suites d'entiers $(i_1,\dots,i_n)$
telles que
\[
0<i_1\le 2i_2\le 4i_3\le \cdots \le 2^{n-1}i_n=2^{n-1}.
\]
En tant qu'espace
vectoriel gradu\'e le module quotient ci-dessus est engendr\'e par les \'el\'ements
$g_n(x_1^{i_1}\cdots x_n^{i_n})$ avec
$(i_1,\dots,i_n)\in\Omega_n$.
En effet, soit un \'el\'ement $g_n(x_1^{a_1}\cdots x_n^{a_n})$  pour
lequel on a $a_{i}> 2 a_{i+1}$. A l'aide de $MP(i)$, il peut
s'\'ecrire comme somme d'\'el\'ements
$g_n(x_1^{a_1}\cdots x_{i-1}^{a_{i-1}}\cdot x_i^{a'_i}\cdot x_{i+1}^{a'_{i+1}}
\cdot x_{i+2}^{a_{i+2}}
\cdots x_k^{a_n})$ avec $a'_i+a'_{i+1}=a_i+a_{i+1}$ et $a'_i<a_i$:
$$x_i^{a_i}x_{i+1}^{a_{i+1}}\equiv \sum_{j=1}^{a_{i+1}}(\binom{a_{i+1}}{j}+
\binom{a_i-a_{i+1}}{j})x_{i}^{a_{i+1}+j}x_{i+1}^{a_i-j} +
 \sum_{j=a_{i+1}+1}^{a_i-a_{i+1}-1}\binom{a_i-a_{i+1}}{j}x_i^{a_{i+1}+j}x_{i+1}^{a_i-j}.
$$
Une
it\'eration \'evidente -tenant compte de $MP(n)$- donne le r\'esultat.

Utilisant la s\'erie de Poincar\'e \ref{poincare-J} de $J(2^n-1)$, 
on observe qu'en tout degr\'e la dimension de l'image de $g_n$ est inf\'erieure ou \'egale \`a celle de $J(2^n-1)$. 
On obtient alors l'isomorphisme souhait\'e.
\end{proof}

\subsection{Surjectivit\'e de $g_n$}
Le reste de la section est consacr\'e \`a la d\'emonstration de la surjectivit\'e de $g_n$.
On proc\`ede comme suit. Soit $V$ un espace vectoriel de dimension $n$. On va
montrer qu'il y a un morphisme surjectif de $H^*(V)$ vers $J(2^n-1)$, puis
on montrera que $$
{\rm Hom}_{\cal U}(H^*(V), J(2^n-1)) \cong H_{2^n-1}(V)$$ est un $\F[{\rm End}(V)]$-module
cyclique de g\'en\'erateur  $g_n$.  La surjectivit\'e de $g_n$ est alors \'evidente.

Pour la premi\`ere partie de l'argument on se sert  de l'action tordue, introduite par
N. Campbell et P. Selick \cite{CS90}, de l'alg\`ebre de Steenrod
sur l'alg\`ebre polynomiale.

On consid\`ere l'alg\`ebre polynomiale $\F[t_0,\ldots,t_{n-1}]$, $t_i$ \'etant de degr\'e $1$.
L'action tordue de l'alg\`ebre de Steenrod sur celle-ci est d\'etermin\'ee par :
$$\Sq^1(t_i)=t_{i-1}^2, \quad 1\le i\le n-1, \quad \Sq^1 (t_0)=t_{n-1}^2$$
et la formule de Cartan. Campbell et Selick montrent alors que, en tant que modules
instables,
les alg\`ebres $\F[t_0,\ldots,t_{n-1}]$ (avec l'action tordue de $\A$) et 
$\F[x_1,\ldots,x_n]$ (munie de l'action classique de $\A$) sont isomorphes.

On introduit une  bigraduation sur $\F[t_0,\ldots,t_{n-1}]$ en imposant que pour
chaque $t_i$  le second degr\'e soit $w(t_i)=2^i$ (comme pour $J_*^*$ plus haut).

Le module instable $\F[t_0,\ldots,t_{n-1}]$ admet alors une d\'ecomposition
en somme directe de $2^n-1$ sous-modules instables, chaque facteur,
soit $\F[t_0,\ldots,t_{n-1}]_i$, \'etant
le sous-module engendr\'e par les mon\^omes dont le second degr\'e est congru modulo $2^n-1$ \`a $i$.
On observe  si $f\in \F[t_0,\ldots,t_{n-1}]$ alors
$$\Sq^1(t_0f)=t_0\Sq^1(f)+t_n^2f, \quad w(t_n^2f)=w(t_0f)+2^n-1$$
donc le sous-espace vectoriel gradu\'e engendr\'e par les mon\^omes dont le second degr\'e est sup\'erieur
\`a une valeur donn\'ee est un sous-module instable.

Pour un \'el\'ement $f\in J^*_* \cong \F[\hat t_i\mid i\geq 0]$ on a
$$\Sq^1(t_0f)=t_0\Sq^1(f).$$

On consid\`ere alors la surjection \'evidente qui envoie $\F[t_0,\ldots,t_{n-1}]$ sur le module de Brown-Gitler $J(2^n-1)\subset J_*^* $ : elle envoie
sur $0$, les mon\^omes de degr\'e sup\'erieur \`a $2^n-1$, et ceux de second degr\'e non congruents \`a $2^n-1$.
On peut voir \'egalement cette application comme \'etant la compos\'ee de l'application d'alg\`ebre de $\F[t_0,\ldots,t_{n-1}]$  vers $J^*_*$ qui envoie $t_i$ vers $\hat t_i$, suivie de la projection sur $J(2^n-1)$.

Si la premi\`ere application n'est pas $\A$-lin\'eaire on v\'erifie facilement que la compos\'ee l'est 
par les arguments donn\'es ci-dessus.

On obtient ainsi un \'epimorphisme de $H^*V$ sur $J(2^n-1)$,
$V$ \'etant un espace vectoriel de dimension $n$.

On \'etudie maintenant le module d'homologie $H_{2^n-1}(V)$.
L'alg\`ebre polynomiale $S^*(V^*)=\F[x_1,\ldots,x_n]$
s'identifie \`a la cohomologie $H^*(V)$ et, de mani\`ere duale, l'alg\`ebre \`a puissances
divis\'ees $\Gamma^*(V)=\Gamma(a_1,\ldots,a_n)$ s'identifie \`a l'homologie $H_*(V)$. Les matrices
\`a coefficients dans $\F$ op\`erent \`a gauche sur $\F[x_1,\ldots,x_n]$, donc par dualit\'e
\`a droite sur $\Gamma(a_1,\ldots,a_n)$, par substitution lin\'eaire des g\'en\'erateurs.

Pour toute suite d'entiers $I=(i_1,\ldots,i_n)$, notons $X^I=x_1^{i_1}\cdots x_n^{i_n}$ et
$A^{(I)}=a_1^{(i_1)}\cdots a_n^{(i_n)}$: les $X^I$ et les $A^{(I)}$ forment respectivement
des bases duales de $\F[x_1,\ldots,x_n]$ et  $\Gamma(a_1,\ldots,a_n)$.
Un multi-indice ou  mon\^ome, $I$, $X^I$ $A^{(I)}$ qui v\'erifie $i_{s-1}\geq 2i_s$ pour $1< s\le n$ sera dit admissible.

\begin{proposition}\label{cyclic} Le $M_n$-module $\Gamma^{2^n-1}(V)$ est engendr\'e par
$a_1^{(2^{n-1})}\cdots a_{n-1}^{(2)} a_n^{(1)}$.  De mani\`ere \'equivalente, pour tout \'el\'ement
$P$ de $S^{2^n-1}(V^*)$, il existe $\sigma\in \m_n(\F)$ tel que l'expression de $\sigma\cdot P$
comme somme de mon\^omes distincts contienne $x_1^{2^{n-1}}\cdots x_{n-1}^2x_n$.
\end{proposition}

La d\'emonstration de l'\'equivalence des deux \'enonc\'es est laiss\'ee au lecteur.

Mettons l'ordre lexicographique sur les mon\^omes de $\F[x_1,\ldots,x_n]$. Pour tout \'el\'ement homog\`ene
non-nul $P\in \F[x_1,\ldots,x_n]$, soit $m(P)$ le plus grand mon\^ome
(par rapport \`a l'ordre lexicographique) qui appara\^it dans $P$.

On aura besoin du lemme suivant :

\begin{lemma}\label{maximal}
Si $m(P)$ n'est pas admissible, alors il existe $\sigma\in \m_n(\F)$ tel que
$m(\sigma\cdot P)> m(P)$.
\end{lemma}

\begin{proof}
Montrons la proposition par r\'ecurrence sur $n$. 
On passera au lemme apr\`es. Le cas $n=1$ est trivial. Supposons que $n>1$
et que l'\'enonc\'e est vrai pour $n-1$. Cette hypoth\`ese implique que pour tout
$Q\in \F[x_2,\ldots,x_n]$ non nul de degr\'e $2^{n-1}-1$, il existe $\tau\in M_n$ tel que
$\tau \cdot Q$ contient $x_2^{2^{n-2}}\cdots x_{n-1}^2x_n$. Ici $\tau$ ne fait intervenir
que les g\'en\'erateurs $x_2,\ldots,x_n$.

Soit $P\in \F[x_1,\ldots,x_n]$ un \'el\'ement quelconque non-nul de degr\'e $2^n-1$. Il faut montrer
que $\sigma \cdot P$ contient $x_1^{2^{n-1}}\cdots x_{n-1}^2x_n$ pour un certain $\sigma\in M_n$.
D'apr\`es le lemme \ref{maximal}, on peut supposer que $m(P)=x_1^{i_1}\cdots x_n^{i_n}$ est un mon\^ome
admissible. 
Si
$m(P)$ est admissible c'est clair. Dans le cas contraire,
en appliquant plusieurs fois le lemme \ref{maximal}, 
on trouve un  $\sigma \in M_n(\F)$  tel que $m(P)$
est admissible. 

Mais alors $i_1\geq 2^{n-1}$.
R\'e\'ecrivons $P$ sous la forme $P=x_1^{2^{n-1}}f+R$, o\`u$f=f(x_1,\ldots,x_n)$ et $R$ ne contient
que des mon\^omes dont la puissance de $x_1$ est inf\'erieure \`a $2^{n-1}$.

Soit $u$ une combinaison lin\'eaire non nulle des g\'en\'erateurs $x_2,\ldots,x_n$. Soit $\sigma_u$
la matrice d\'efinie par la substitution $x_1:=x_1+u$. Posons $Q=f(u,x_2,\ldots,x_n)$. Le polyn\^ome, en $x_2,\ldots,x_n$, coefficient de $x_1^{2^{n-1}}$ dans $$\sigma_u\cdot P=P(x_1+u,x_2,\ldots,x_n)
\in \F[x_2,\ldots,x_n][x_1]$$ est $Q$. On suppose d'abord que  $Q\not =0$.
Or, par r\'ecurrence il existe $\tau\in M_n$ tel que
$\tau\cdot Q$ contient $x_2^{2^{n-2}}\cdots x_{n-1}^2x_n$. Comme
$\tau$ ne fait pas intervenir $x_1$, il suit que $\tau\sigma_u\cdot P$ contient
$x_1^{2^{n-1}}\cdots x_{n-1}^2x_n$.

Si  $f(u,x_2,\ldots,x_n)=0$ pour toute combinaison lin\'eaire non nulle $u$,
alors  $f(x_1,\ldots,x_n)$ est divisible par $x_1+u$ quelque soit $u$.
Comme $f$ est de degr\'e $2^{n-1}-1$, il suit que $f=\prod_u(x_1+u)$ et
$x_2^{2^{n-2}}\cdots x_{n-1}^2x_n$ appara\^it dans $f$. D'o\`u $P$ lui-même contient
$x_1^{2^{n-1}}\cdots x_{n-1}^2x_n$,
la proposition est d\'emontr\'ee.
\end{proof}

\proofof{\ref{maximal}}
On met l'ordre lexicographique sur les mon\^omes de $\F[x_1,\cdots,x_n]$.
Pout tout \'el\'ement  non-nul $P \in \F[x_1,\cdots,x_n]$,
soit $m(P)$ le plus grand mon\^ome (par rapport \`a l'ordre lexicographique)
qui appara\^it dans $P$.

 Soit $m(P) = x^{i_1}\cdots x^{i_n}_n$. Comme $m(P)$ n'est pas admissible, il
existe $1 < s \leq n$ tel que $2i_{s} > i_{s-1}$. En regroupant les mon\^omes de $P$, on le
r\'e\'ecrit sous la forme $
P = x^{i_1}_1
\cdots x^{i_{s-2}}
_{s-2}Qx^{i_{s+1}}_{s+1}\cdots
x^{i_{n}}_n+ R$,
de mani\`ere que $m(R) < m(P)$ et que $Q = Q(x_{s-1}, x_s) \in \F[x_{s-1}, x_s]$ est un
polyn\^ome  de degr\'e $i_{s-1} + i_s$ qui v\'erifie $m(Q) = x^{i_{s-1}}_{s-1} x^{i_s}_s$ .

Pour tout $\sigma \in  {\rm GL}_n$ qui correspond \`a une  substitution qui
ne fait intervenir que  $x_{s-1}, x_s$, les mon\^omes de $\sigma \cdot R $ sont
diff\'erents des mon\^omes de
$
x^{i_1}_1
\cdots x^{i_{s-2}}
_{s-2} (\sigma \cdot Q)x^{i_{s+1}}_{s+1}\cdots
x^{i_{n}}_n$. Si on a $m(\sigma  \cdot  x^{i_1}_1
\cdots x^{i_{s-2}}
_{s-2} (\sigma \cdot Q)x^{i_{s+1}}_{s+1}\cdots
x^{i_{n}}_n )> m(x^{i_1}_1
\cdots x^{i_{s-2}}
_{s-2} (\cdot Q)x^{i_{s+1}}_{s+1}\cdots
x^{i_{n}}_n)$   on aura forc\'ement  $m(\sigma \cdot P) > m(P)$.

On pose $Q = x^q_{s-1}x^q_s(x_{s-1} + x_s)^qQ'$ avec $Q' \in  \F[x_{s-1}, x_s]$ et l'entier $q$ le
plus grand possible.

Dans le premier cas, $Q'$ contient un mon\^ome $x^i_{s-1}$. Celui-ci
est forc\'ement $m(Q')$, d'o\`u $m(Q) = x^{2q+i}_{s-1}x_s^q$, ce qui est absurde puisque
$m(Q) = x_{s-1}^{i_{s-1}} x^{i_ s}_s$ et $2i_s > i_{s-1}$.

Dans le deuxi\`eme cas, $Q'$ contient un mon\^ome $x^i_{s}$. Soit $\sigma \in {\rm GL}_n$ la matrice qui
correspond \`a la substitution qui \'echange $x_{s-1}$ et  $x_s$. Alors
$$\sigma \cdot Q = Q(x_s, x_{s-1}) =
x^q_{s-1}x_q^s(x_{s-1} + x_s)^qQ'(x_s, x_{s-1)}.
$$
 Il est clair que $x^i_{s-1} = m(\sigma \cdot Q')$, d'o\`u $ m(\sigma \cdot Q')=x^{2q+i}_{s-1}x_s^q$. Comme $3q+i = i_{s-1}+i_s$, on d\'eduit ais\'ement que $m(Q) > x^{i_{s-1}}_{s-1} x^{i_ s}_s$,
d'o\`u $ m(\sigma \cdot P) > m(P)$.

Dans le troisi\`eme cas, $Q'$ est divisible par $x_{s-1}x_s$. Soit $\tau \in {\rm GL}_n$ la matrice qui
correspond \`a la substitution qui transforme $x_s$ en  $x_{s-1}+x_s$. Alors dans
$$\tau \cdot Q' = Q'(x_{s-1},x_{s-1} + x_s)$$
le terme qui ne comporte pas $x_s$  est \'egal \`a
$Q'(x_{s-1}, x_ {s-1})$. Puisque $Q'$ n'est pas divisible par $x_{s-1}+x_s$ par maximalit\'e de $q$, ce terme est non nul
et  \'egal \`a $x^i_{s-1}$  pour un certain $i > 0$. Il suit que $m(\tau \cdot Q) = x^{2q+i}_{s-1}x_s^q$.
Comme $3q + i = i_{s-1} + i_s$, on a $m(\tau \cdot Q) > x^{i_{s-1}}_{s-1}
 x^{i_s}_s$, d'o\`u $m(\tau \cdot P)> m(P)$.
\qed

%******************************************************

\section{Applications homologiques}\label{application-homology}

Dans cette section on  d\'emontre diverses cons\'equences homologiques  du th\'eor\`eme \ref{NST1}.

On commence par
\begin{theorem} La r\'esolution injective de $L'_n$ donn\'ee par \ref{NST1}
est minimale.
\end{theorem}

Cela r\'esulte de ce que les modules instables $L_h \otimes J(k)$ sont 
ind\'ecomposables et deux \`a deux distincts \cite{LS89}. 

\begin{corollary} \label{mu}
L'\'el\'ement $\mu_n \in {\rm Ext}^n_{\cal U}(J(2^n-1),L'_n)$ 
d\'etermin\'e par cette r\'esolution injective est non nul.
\end{corollary}

Le corollaire r\'esulte de ce que $J(2^n-1)$ est localement fini, 
alors que le plus grand sous-module localement fini de $L_1 \otimes J(2^{n-1}-1)$ est trivial. 
En effet ceci est cons\'equence de ce que $L_1 \otimes J(2^{n-1}-1)$ a une filtration finie dont les quotients 
sont des suspensions de $L_1$, voir aussi \cite[Chapter 6]{Sch94}.

\begin{corollary} Soit $M$ un module instable, $n>0$, 
\begin{enumerate}
\item ${\rm Ext}^s_{\cal U}(M,L'_n)=\{0\}$ si $s >n$;
\item ${\rm Ext}^s_{{\mathcal U}}(M,L'_n)=\{0\}$ si $s \not = n$ et $M$ est localement fini;
\item ${\rm Ext}^s_{{\mathcal U}}(\Sigma^t M,L'_n)=\{0\}$  si $t > 2^s-1$.
\end{enumerate}
\end{corollary}

\begin{proof}
La premi\`ere propri\'et\'e est claire, la seconde   vient de ce que ${\rm Hom}_{\cal U}(M,L_j \otimes J(k)) $
est nul si $j>0$ car $M$ est localement fini et $L_j \otimes J(k)$ a une partie
localement finie triviale si $j>0$. La troisi\`eme r\'esulte de ce qu'il n'y a pas d'applications non nulles 
de $\Sigma^h M$ dans $L_k$ si $h>0$. 
\end{proof}

\subsection{Groupes d'extension $\Ext_\A^{s,t}(\Z/2,L'_n)$, $n\ge 2$}

On consid\`ere ensuite le calcul des groupes d'extension $\Ext_{\A}^{s,t}(\Z/2,L_n')$. 
A cet effet, on se sert des travaux de J. Lannes et S. Zarati \cite{LZ87} sur les foncteurs d\'eriv\'es 
de la d\'estabilisation. 

On d\'esigne par $\Mc$ la cat\'egorie dont les objets sont les $\A$-modules gradu\'es et 
dont les morphismes sont les applications $\A$-lin\'eaires de degr\'e z\'ero. La cat\'egorie $\U$ des 
$\A$-modules instables est alors une sous cat\'egorie pleine de $\Mc$. 
On note $D\colon \Mc\rightarrow \U$ et on appelle {\it foncteur de d\'estabilisation} 
l'adjoint \`a gauche du foncteur oubli $\U \rightarrow \Mc$.  
Le foncteur $D$ est exact \`a droite et on note 
$D_s\colon \Mc\rightarrow \U$, $s\ge 0$, ses foncteurs d\'eriv\'es. 

Pour tout module instable $M$ et tout $s\ge 0$, Lannes et Zarati 
ont explicit\'e le module $D_s\Sigma^{1-s}M$ en termes de la construction de Singer de $M$. 
En particulier, 
si $D(s)$ et $\omega_s$ d\'esigne respectivement l'alg\`ebre de Dickson et 
l'invariant de Dickson sup\'erieur de $\gl_s$ sur $\F[x_1,\cdots,x_s]$, 
on a
\begin{proposition}[\cite{LZ87}]\label{LZformula}
$D_s\Sigma^t(\Z/2)\cong \Sigma^{s+t}D(s)\omega_s^{s+t-1}$ 
si $s\ge 0$ et $s+t\ge 1$. 
\end{proposition}
\begin{proof} 
Soit $M$ un module instable. D'après Lannes et Zarati \cite{LZ87}, on a 
$D_s\Sigma^{1-s} M\cong\Sigma R_s M$, $R_s(\Sigma M)\cong\Sigma \omega_s R_s M$ et $R_s(\Z/2)\cong D(s)$. Ici 
$R_s(M)$ d\'esigne un certain sous-$\A$-module de $H^*B(\Z/2)^s\otimes M$ d\'ependant 
fontoriellement de $M$. 
On en deduit que $$D_s\Sigma^t(\Z/2)=D_s\Sigma^{1-s} (\Sigma^{s+t-1}\Z/2)\cong
\Sigma R_s(\Sigma^{s+t-1}\Z/2)\cong\Sigma^{s+t}\omega_s^{s+t-1}R_s(\Z/2)\cong\Sigma^{s+t}D(s)\omega_s^{s+t-1}.$$ 
La proposition est démontrée.
\end{proof}

Soit $N$ un module instable. La suite spectrale de Grothendieck associ\'ee \`a la compos\'ee des foncteurs 
$$\xymatrix{\Mc\ar[rr]^D && \U\ar[rr]^-{\Hom_\U(-,N)} && \text{$\F$-espaces vectoriels}}$$
est de la forme :
$$E_2^{p,q}:=\Ext_\U^p(D_q(-),N)\Rightarrow \Ext_\Mc^{p+q}(-,N).$$

\begin{figure}[ht]\label{E2}
\centering
\begin{pspicture}(-8,-2)(6,7)
\psline{->}(0,0)(0,7)
\psline{->}(-8,0)(6,0)
\pspolygon[fillcolor=brown,fillstyle=solid](-8,0)(-1,0)(-1,7)(-7,7)
\uput[d](6,0){$t-s$}
\uput[r](0,7){$s$}
\qdisk(-5,2){2pt}
\uput[ur](-5,2){$\Z/2$}
\uput[ur](-2.5,3){$\Z/2$}
\uput[d](0,0){$0$}
\uput[d]{90}(-1,0){$-n-2^{n-2}$}
\uput[d]{90}(-5,0){$-n-2^n+1$}
\uput[d]{90}(-2.5,0){$-n-2^{n-1}$}%
\uput[r](0,3){$n+1$}
\uput[r](0,2){$n$}
\rput(-6.5,3){\psframebox*[framearc=0.4]{\Huge $\mathfrak A_n$}}
\psline{->}(-4.7,1)(-5,2)\psline{->}(-5,2)(-5.3,3)
\psline{->}(-2.2,2)(-2.5,3)\psline{->}(-2.5,3)(-2.8,4)
\psline[linecolor=red,linestyle=dotted](-5.3,0)(-5.3,7)
\psline[linecolor=red,linestyle=dotted](-4.7,0)(-4.7,7)
\psline[linecolor=red,linestyle=dotted](-5,0)(-5,7)
\psline[linecolor=red,linestyle=dotted](0,2)(-5,2)
\psline[linecolor=red,linestyle=dotted](-2.5,0)(-2.5,3)(0,3)
\qdisk(0,2){2pt}
\qdisk(0,3){2pt}
\qdisk(-1,0){2pt}
\qdisk(-5,0){2pt}
\qdisk(-2.5,0){2pt}
\qdisk(-2.5,3){2pt}
\end{pspicture}
\caption{${\Ext}_{\A}^{s,t}(\Z/2,L'_n)$, $n\geq 2$}
\end{figure}

Le th\'eor\`eme \ref{NST2} est la combinaison des trois propositions suivantes.

\begin{proposition} On a $\Ext_\A^{n,1-2^n}(\Z/2,L'_n)=\Z/2$.
\end{proposition}

\begin{proof}
On a une suite spectrale de Grothendieck:
$$\Ext_\U^p(D_q\Sigma^{2^n-1}\Z/2,L_n')\Rightarrow \Ext_{\Mc}^n(\Sigma^{2^n-1}\Z/2,L_n').$$
D'après \ref{LZformula},
$$D_q\Sigma^{2^n-1}\Z/2= \Sigma^{q+2^n-1}D(q)\omega_q^{q+2^n-2}.$$
Soit $\widetilde \Sigma \colon \U\rightarrow \U$ le foncteur adjoint \`a droite du foncteur 
de suspension $\Sigma \colon \U\rightarrow \U$. On rappelle que 
$\widetilde \Sigma(L\otimes J(k))=L\otimes J(k-1)$ 
si $L$ est un module instable r\'eduit \cite{Sch94}. On adopte la convention que $J(k)=0$ si $k$ est n\'egatif.
On déduit alors de ce qui pr\'ec\`edent et de la résolution injective de $L'_n$ que  
\begin{eqnarray*}
\Ext_\U^p(D_q\Sigma^{2^n-1}\Z/2,L_n')=
\begin{cases}0 &\text{si $p\not =n$,}\\
\Hom_\U(D_q\Sigma^{2^n-1}\Z/2,J(2^n-1)) & \text{si $p=n$,}
\end{cases} & =& \begin{cases}0 &\text{si $(p,q)\not =(n,0)$,}\\
\Z/2 & \text{si $(p,q)=(n,0)$.}
\end{cases}
\end{eqnarray*}
Il suit que la suite spectrale d\'eg\'en\`ere et l'on obtient 
$$\Ext_{\Mc}^n(\Sigma^{2^n-1}\Z/2,L_n')\cong \Ext_\U^n(\Sigma^{2^n-1}\Z/2,L_n')\cong\Z/2.$$
La proposition est d\'emontr\'ee.
\end{proof}

\begin{remark} \label{nu}
On d\'eduit de la d\'emonstration ci-dessus que 
l'inclusion $\Sigma^{2^n-1}\Z/2 \hookrightarrow J(2^n-1)$ induit un isomorphisme
$\Ext_\U^n(J(2^n-1),L'_n)\cong \Ext_\Mc^n(\Sigma^{2^n-1}\Z/2,L'_n)\cong \Z/2$ pour tout $n\ge 1$. 
On peut alors d\'efinir une classe non-nulle $\nu_n\in \Ext_\Mc^n(\Sigma^{2^n-1}\Z/2,L'_n)$ 
comme \'etant l'image de $\mu_n\in \Ext_\U^n(J(2^n-1),L'_n)$ par cet isomorphisme. 
Une suite exacte des $\A$-modules qui repr\'esente $\nu_n$ 
sera donn\'ee dans \ref{cofibration-taka}.
\end{remark}

\begin{proposition}On a $\Ext_\A^{n+1,1-2^{n-1}}(\Z/2,L'_n)=\Z/2$.
\end{proposition}
\begin{proof}
On a une suite spectrale de Grothendieck :
$$\Ext_\U^p(D_q\Sigma^{2^{n-1}-1}\Z/2,L_n')\Rightarrow \Ext_{\Mc}^{n+1}(\Sigma^{2^{n-1}-1}\Z/2,L_n').$$
D'apr\`es \ref{LZformula},
$$D_q\Sigma^{2^{n-1}-1}\Z/2= \Sigma^{q+2^{n-1}-1}D(q)\omega_q^{q+2^{n-1}-2}.$$
\begin{itemize}
\item Si $q=0$ alors, pour tout $0\le i\le n$,
$$\Hom_\U(\Sigma^{2^{n-1}-1}\Z/2,L_{n-i}\otimes J(2^i-1))=0.$$

\item Si $q=1$ alors $D_1\Sigma^{2^{n-1}-1}\Z/2= \Sigma^{2^{n-1}}D(1)\omega_1^{2^{n-1}-1}$, donc
$$
\Ext_\U^p(D_1\Sigma^{2^{n-1}-1}\Z/2,L_n')=
\begin{cases}0 &\text{si $p\not =n$,}\\
\Z/2 & \text{si $p=n$.}
\end{cases}
$$

\item Si $q\geq 2$ alors la connectivit\'e de 
$D_q\Sigma^{2^{n-1}-1}\Z/2=\Sigma^{q+2^{n-1}-1}D(q)\omega_q^{q+2^{n-1}-2}$ est
$$2^q(q+2^{n-1}-2)+1>2^n-1.$$
\end{itemize}
On en d\'eduit que
$$
\Ext_\U^p(D_q\Sigma^{2^{n-1}-1}\Z/2,L_n')=
\begin{cases}
0 &\text{si $(p,q)\not =(n,1)$,}\\
\Z/2 & \text{si $(p,q)=(n,1)$.}
\end{cases}
$$
D'o\`u $\Ext_{\Mc}^{n+1}(\Sigma^{2^{n-1}-1}\Z/2,L_n')\cong \Ext_\U^n(D_1\Sigma^{2^{n-1}-1}\Z/2,L_n')\cong\Z/2$. 
\end{proof}

\begin{proposition} Supposons $(s,t)\not \in \{(n,1-2^n),(n+1,1-2^{n-1})\}$.
Alors le groupe $\Ext_\A^{s,t}(\Z/2,L'_n)$ est nul d\`es que $s-t\ge 2^{n-2}+n$.
\end{proposition}

\begin{proof}
On a une suite spectrale de Grothendieck :
$$\Ext_\U^p(D_{s-p}\Sigma^{-t}\Z/2,L_n')\Rightarrow \Ext_{\Mc}^{s}(\Sigma^{-t}\Z/2,L_n').$$
D'apr\`es \ref{LZformula},
$$D_{s-p}\Sigma^{-t}\Z/2= \Sigma^{s-p-t}D(s-p)\omega_{s-p}^{s-p-t-1}.$$
On va montrer que $\Ext_\U^p(D_{s-p}\Sigma^{-t}\Z/2,L_n')=0$ pour $0\le p\le n$.
\begin{enumerate}
\item Pour $0\le p\le n-2$, on a $s-p-t> 2^{n-2} -1$. 
On en d\'eduit que le groupe $\Ext_\U^p(D_{s-p}\Sigma^{-t}\Z/2,L_n')$ est nul pour tout $0\le p\le n-2$. 

\item Pour $p=n-1$, le groupe $\Ext_\U^{n-1}(D_{s-n+1}\Sigma^{-t}\Z/2,L_n')$ est le 
noyau du morphisme 
$$\Hom_\U\big(D_{s-n+1}\Sigma^{-t}\Z/2,L_1\otimes J(2^{n-1}-1)\big)
\xrightarrow{f} \Hom_\U\big(D_{s-n+1}\Sigma^{-t}\Z/2,J(2^n-1)\big)$$
qui est induit par le morphisme $f_{n-1,n}\colon L_1\otimes J(2^{n-1}-1)\rightarrow J(2^n-1)$. 
Utilisant le fonteur $\widetilde\Sigma$, le morphisme $f$ se r\'e\'ecrit comme suit :
$$\Hom_\U\big(D(s-n+1)\omega_{s-n+1}^{s-t-n},L_1\otimes J(2^{n-1}-2+n+t-s)\big)
\xrightarrow{f'} \Hom_\U\big(D(s-n+1)\omega_{s-n+1}^{s-t-n},J(2^n-2+n+t-s)\big)$$
qui est induit par le morphisme compos\'e 
$$L_1\otimes J(2^{n-1}-2+n+t-s)\xrightarrow{\pi_{2^{n-1}}\otimes id} J(2^{n-1})\otimes J(2^{n-1}-2+n+t-s) \xrightarrow{\mu} J(2^n-2+n+t-s).$$

\begin{enumerate}
	\item Si $s=n-1$, le morphisme $f'$ devient 
		$$\Hom_\U(\F, L_1\otimes J(2^{n-1}-1+t))\rightarrow \Hom_\U(\F,J(2^n-1+t))$$
		avec $t\le -2^{n-2}-1$. Le domaine de ce morphisme est trivial.

	\item Si $s=n$, le morphisme $f'$ devient 
		$$\Hom_\U\big(D(1)\omega_{1}^{-t},L_1\otimes J(2^{n-1}-2+t)\big)
		\rightarrow \Hom_\U\big(D(1)\omega_{1}^{-t},J(2^n-2+t)\big)$$
	avec $t\le -2^{n-2}$. Si $t<2-2^{n-1}$, le domaine du morphisme est nul. 
	On va montrer que ce morphisme est un isomorphisme si 
	$2-2^{n-1}\le t\le -2^{n-2}$. En effet, on a 
		\begin{eqnarray*}
			\Hom_\U\big(D(1)\omega_{1}^{-t},L_1\otimes J(2^{n-1}-2+t)\big)& \cong & 
			\Hom_\U\big(\widetilde T L(1), J(2^{n-1}-2+t)\big)\\
			&\cong&\Hom_\U\big(\F[x], J(2^{n-1}-2+t)\big)\\
			&\cong& \Z/2.
		\end{eqnarray*} 
	Ici $\widetilde T$ 
	est le foncteur de Lannes qui est adjoint \`a gauche du foncteur $- \otimes L_1 \colon \U\rightarrow \U$.
	L'\'el\'ement non trivial, not\'e $\alpha$, de $\Hom_\U\big(D(1)\omega_{1}^{-t},L_1\otimes J(2^{n-1}-2+t)\big)$ 
	s'\'ecrit comme \'etant le compos\'e :
	$$D(1)\omega_{1}^{-t} \xrightarrow{\iota }\F[x,y] \rightarrow 
	\F[x]\otimes \F[y] \rightarrow L_1\otimes J(2^{n-1}-2+t),$$
	où $\iota(x^k)=(x+y)^k$. Or, le coefficient de $x^{2^{n-1}}y^{2^{n-1}-2+t}$ 
	dans le d\'eveloppement de $(x+y)^{2^n-2+t}$ est $\binom{2^n-2+t}{2^{n-1}}$ qui est non nul puisque 
	$2^{n-1} \le 2^n-2+t \le 2^{n-1}+2^{n-2}-2$. Le compos\'e 
	$$D(1)\omega_{1}^{-t} \xrightarrow{\alpha} L_1\otimes J(2^{n-1}-2+t) \xrightarrow{\pi_{2^{n-1}}\otimes id} 
	J(2^{n-1})\otimes J(2^{n-1}-2+t) \rightarrow J(2^n-2+t)$$ 
	est ainsi non nul car il envoie $x^{2^n-2+t}$ sur $\iota_{2^n-2+t}$.

	\item Si $s\ge n+1$, la connectivit\'e de $\widetilde T(D(s-n+1)\omega_{s-n+1}^{s-t-n})$ 
	est $(2^{s-n}-1)(s-t-n)$ d'apr\`es la proposition \ref{connectivity} ci-dessous. 
	On en déduit que 
	le groupe $\Hom_\U\big(D(s-n+1)\omega_{s-n+1}^{s-t-n},L_1\otimes J(2^{n-1}-2+n+t-s)\big)$ 
	est nul car 
	$$(2^{s-n}-1)(s-t-n)=2^{s-n}(s-t-n)+n+t-s>2^{n-1}-2+n+t-s.$$
	\end{enumerate}

\item Pour $p=n$, comme le cas pr\'ec\'edent, 
le groupe $\Ext_\U^n(D_{s-n}\Sigma^{-t}\Z/2,L_n')$ est le conoyau du morphisme 
$$\Hom_\U\big(D(s-n)\omega_{s-n}^{s-t-n-1}, L_1\otimes J(2^{n-1}+n-1+t-s)\big)\xrightarrow{f''}
\Hom_\U\big(D(s-n)\omega_{s-n}^{s-t-n-1}, J(2^n+n-1+t-s)\big).$$
	\begin{enumerate}
		\item	Si $s=n$, le morphisme $f''$ devient 
						$$\Hom_\U(\F,L_1\otimes J(2^{n-1}-1+t)\rightarrow \Hom_\U(\F,J(2^n-1+t))$$
				avec $t\le 1-2^{n-2}$. 
		La source de ce morphisme est triviale car $t\not= 1-2^n$.
		\item   Si $s=n+1$, le morphisme $f''$ devient 
						$$\Hom_\U\big(D(1)\omega_{1}^{-t}, L_1\otimes J(2^{n-1}-2+t)\big)\rightarrow
				\Hom_\U\big(D(1)\omega_{1}^{-t}, J(2^n-2+t)\big)$$
				avec $t\le 1-2^{n-2}$. Le groupe $\Hom_\U\big(D(1)\omega_{1}^{-t}, J(2^n-2+t)\big)$ 
				est non trivial si et seulement si $t\ge 1-2^{n-1}$. Comme $t\not = 1-2^{n-1}$, le morphisme 
				consid\'er\'e avec $1-2^{n-1}<t\le 1-2^{n-2}$ est un isomorphisme comme le cas $2(b)$ ci-dessus.
		\item 	Si $s=n+2$, le morphisme $f''$ devient 
			$$\Hom_\U\big(D(2)\omega_{2}^{-t+1}, L_1\otimes J(2^{n-1}-3+t)\big)\rightarrow
			\Hom_\U\big(D(2)\omega_{2}^{-t+1}, J(2^n-3+t)\big)$$
				avec $t\le 2-2^{n-2}$. Le module $D(2)\omega_{2}^{-t+1}$ est alors non nul 
				en degr\'e $2^n-3+t$ si et seulement si $t=2-2^{n-2}$. Le morphisme devient  
				$$\Hom_\U\big(D(2)\omega_{2}^{2^{n-2}-1}, L_1\otimes J(2^{n-2}-1)\big)\rightarrow
			\Hom_\U\big(D(2)\omega_{2}^{2^{n-2}-1}, J(3.2^{n-2}-1)\big)\cong \Z/2.$$
				Ce morphisme est surjectif car l'application 
				$D(2)\omega_{2}^{2^{n-2}-1}\hookrightarrow L_1\otimes L_1
				\xrightarrow{id\otimes \pi_{2^{n-2}-1}} L_1\otimes J(2^{n-2}-1)$ envoie 
				la classe $(x^2+xy+y^2)\omega_{2}^{2^{n-2}-1}$ sur la classe 
				$x^{2^{n-1}}\otimes \iota_{2^{n-2}-1}$.
		\item 	Si $s\ge n+3$, la source du morphisme $f''$ est triviale car $D(s-n)\omega_{s-n}^{s-t-n-1}$ est 
trival en degr\'e $2^n+n-1+t-s$. En effet, supposons qu'il existe un \'el\'ement de $D(s-n)$ de degr\'e $d$ 
tel que $$(2^{s-n}-1)(s-t-n-1)+d=2^n+n-1+t-s.$$ Il suit 
$2^n-2=2^{s-n}(s-t-n-1)+d\ge 8(2^{n-2}-1)+d \Rightarrow 6\ge 2^n+d \Rightarrow d=2$. Or, l'alg\`ebre de Dickson $D(s-n)$ 
est trivial en degr\'e $2$ si $s-n\ge 3$.
	\end{enumerate}
\end{enumerate}
On a ainsi démontré que $\Ext_\U^p(D_{s-p}\Sigma^{-t}\Z/2,L_n')=0$ pour $0\le p\le n$.
La proposition suit.
\end{proof}

\subsection{S\'erie de Poincar\'e de $\widetilde T\big(D(n)\omega_n^i\big)$}
Soit $\Phi$ le foncteur double \cite{Sch94} de la cat\'egrorie $\U$ des modules instables. 
Rappelons que l'alg\`ebre des invariants $\F[x_1,\cdots,x_n]^{\gl_n}$ est l'alg\`ebre 
de Dickson $D(n)=\F[Q_{n,0},\cdots, Q_{n,n-1}]$. Le 
lemme suivant est facile \`a vérifier en utilisant le fait que, 
par la projection canonique $\F[x_1,\cdots,x_n]\twoheadrightarrow \F[x_1,\cdots,x_{n-1}]$, 
l'invariant de Dickson $Q_{n,i}$ s'envoie sur $Q_{n-1,i-1}^2$ si $i>0$ et sur $0$ si $i=0$. 

\begin{lemma} \label{dickson-exact}
Soient $n,i\ge 1$. Il existe une suite exacte courte des modules instables :
$$0\rightarrow D(n)\omega_n^i\rightarrow D(n)\omega_n^{i-1} 
\rightarrow \Sigma^{i-1}\Phi D(n-1)\omega_{n-1}^{i-1}\rightarrow 0.$$
\end{lemma}

On note $P_{n,i}(t)$ la s\'erie de Poincar\'e de $\widetilde T\big(D(n)\omega_n^i\big)$. 
Il r\'esulte du lemme 
\ref{dickson-exact} et de l'exactitude du foncteur $\widetilde T$ que
\begin{equation}\label{Pni}
P_{n,i}(t)=P_{n,i-1}(t)-t^{i-1}P_{n-1,i-1}(t^2).
\end{equation}
\begin{proposition}\label{connectivity}
 On a $P_{n,i}(t)=t^{(2^{n-1}-1)i}P_{n,0}(t)$.
\end{proposition} 
\begin{proof}
Soit $G_{n-1}$ le sous-groupe de $\gl_n$ des matrices de la forme 
$$\begin{pmatrix}* & \cdots & * & * \\ \vdots & \ddots & \vdots & \vdots\\
 * & \cdots & * & * \\
0 & \cdots & 0 &1 \end{pmatrix}$$
La th\'eorie de Lannes \cite{Sch94} donne 
$$\widetilde T\big(D(n)\big )= \F[x_1,\cdots,x_n]^{G_{n-1}}.$$
Or, c'est un résultat classique que l'alg\`ebre des invariants $\F[x_1,\cdots,x_n]^{G_{n-1}}$ 
est une alg\`ebre polynomiale engendr\'ee par les g\'en\'erateurs $Q_{n-1,0},\cdots, Q_{n-1,n-2},V_n$ 
dont les degr\'es sont $2^{n-1}-1,\cdots,2^{n-1}-2^{n-2},2^{n-1}$ respectivement. Il suit 
$$P_{n,0}(t)=\frac{1}{(1-t^{2^{n-1}-1})\cdots (1-t^{2^{n-1}-2^{n-2}})(1-t^{2^{n-1}})}.$$
On en déduit que 
\begin{equation}\label{Pdouble}
P_{n-1,0}(t^2)=(1-t^{2^{n-1}-1})P_{n,0}(t).
\end{equation}
La proposition est maintenant facile \`a d\'emontrer par r\'ecurrence double sur $(i,n)$. 
On a rien \`a faire pour $i=0$. Pour $n=1$, on a $\widetilde T(D(1)\omega_1^i)=\widetilde T(D(1))$ 
car le module quotient $D(1)/D(1)\omega_1^i$ est fini.

Supposons que la formule $P_{n',i'}(t)=t^{(2^{n'-1}-1)i'}P_{n',0}(t)$ soit v\'erifi\'ee pour 
tout $(n',i')<(n,i)$. 
On a
\begin{align*}
P_{n,i}(t)&=P_{n,i-1}(t)-t^{i-1}P_{n-1,i-1}(t^2) && \text{(d'apr\`es \ref{Pni})}
\\
&= t^{(2^{n-1}-1)(i-1)}P_{n,0}(t)-t^{i-1} t^{2(2^{n-2}-1)(i-1)}P_{n-1,0}(t^2) && 
\text{(d'apr\`es l'hypoth\`ese de r\'ecurrence)} \\
&=t^{(2^{n-1}-1)(i-1)}P_{n,0}(t)-t^{(2^{n-1}-1)(i-1)}(1-t^{2^{n-1}-1})P_{n,0}(t) 
&& \text{(d'apr\`es \ref{Pdouble})}\\
&= t^{(2^{n-1}-1)i}P_{n,0}(t). && 
\end{align*}
La proposition est d\'emontr\'ee.
\end{proof}

\begin{remark}
La proposition sugg\`ere que l'on ait un isomorphisme de modules instables:
$$\widetilde T\big(D(n)\omega_n^i\big)\cong\widetilde T\big(D(n)\big)\cdot \omega_{n-1}^i.$$
\end{remark}

\section{Applications homotopiques}\label{application-homotopy}

Dans cette section, tous les espaces et spectres sont $2$-compl\'et\'es.

\subsection{Les cofibrations de Takayasu}
\label{cofibration-taka}
Soit $V_n=(\Z/2)^n$. On consid\`ere la repr\'esentation r\'eelle r\'eguli\`ere r\'eduite de $V_n$:
 $$\reg_n\colon V_n \rightarrow Aut(\R^{2^n-1})$$
et la somme $\reg_n^{\oplus k}=\reg_n+\cdots+\reg_n$, $k$ fois. Takayasu \cite{Tak99} consid\`ere 
$\reg_n^{\oplus k}$ pour tout entier $k$. On se restreint aux cas o\`u $k\ge 0$. On d\'esigne par 
$BV_n^{\reg_n^{\oplus k}}$ l'espace de Thom associ\'e \`a la repr\'esentation $\reg_n^{\oplus k}$, {\it i.e.} l'espace de 
Thom du fibr\'e vectoriel $EV_n\times_{V_n}\R^{k(2^n-1)}\rightarrow BV_n$.

L'action de $\gl_n$ sur $BV_n$ induit une action sur l'espace de Thom 
$BV_n^{\reg_n^{\oplus k}}$. 
L'idempotent de Steinberg $e_n$ d\'efinie alors une application stable sur 
$BV_n^{\reg_n^{\oplus k}}$. A la suite de Mitchell et Priddy \cite{MP83}, en prenant 
le t\'el\'escope de cette application,  on obtient un facteur 
stable de $BV_n^{\reg_n^{\oplus k}}$ que l'on note $e_n\cdot BV_n^{\reg_n^{\oplus k}}$. 
On adopte la notation de Takayasu en posant $\M(n)_k=e_n\cdot BV_n^{\reg_n^{\oplus k}}$. 
On renvoie \`a \cite{Tak99} pour une construction explicite de $\M(n)_k$,. 

On obtient en particulier que $\M(n)_0=\M(n)$, $\M(n)_1=\L(n)$ et  $\M(n)_2=\L'(n)$. 
La cohomologie de $\M(n)_k$ est d\'etermin\'ee en utilisant l'isomorphisme de Thom :
$$H^*\M(n)_k=\omega_n^ke_n\cdot \F[x_1,\cdots,x_n]=\omega_n^{k-1}L_n.$$

\begin{theorem}[Takayasu \cite{Tak99}] Pour $k\ge 0$, il existe une suite de cofibration 
$$\Sigma^k\M(n-1)_{2k+1}\xrightarrow{i_{n,k}} \M(n)_k\xrightarrow{j_{n,k}} \M(n)_{k+1}.$$
\end{theorem}

On obtient en particulier une suite exacte courte en cohomologie : 
$$0\rightarrow \omega_n^{k} L_n \rightarrow \omega_n^{k-1} L_n \rightarrow \Sigma^k\omega_{n-1}^{2k} L_{n-1}\rightarrow 0.$$
On combine les suites de cofibration de Takayasu pour obtenir la suivante :
{\small
\begin{equation}\label{T}\tag{{\bf T}}
\Sigma^{2^n-1}\M(0)_{2^n}\rightarrow \cdots 
\rightarrow \Sigma^{2^k-1} \M(n-k)_{2^k}\xrightarrow{d_{k,n}} \Sigma^{2^{k-1}-1} \M(n-k+1)_{2^{k-1}}\rightarrow
\cdots \rightarrow \Sigma \M(n-1)_2\rightarrow \M(n)_1\rightarrow \M(n)_2.
\end{equation}
}
Ici $d_{k,n}$ se factorise comme suit:
{\small
$$
\xymatrix
{\Sigma^{2^k-1}\M(n-k)_{2^k}\ar[rr]^{d_{k,n}}\ar[dr]^{\Sigma^{2^k-1}j_{n-k,2^k}} & 	&	\Sigma^{2^{k-1}-1} \M(n-k+1)_{2^{k-1}}\\
& \Sigma^{2^k-1}\M(n-k)_{2^k+1}
\simeq  \Sigma^{2^{k-1}-1}\big (\Sigma^{2^{k-1}}\M(n-k)_{2^k+1}\big )\ar[ur]^{\Sigma^{2^{k-1}-1}j_{n-k+1,2^{k-1}}}.
}$$
}
La suite \ref{T} induit en cohomologie une suite exacte de $\A$-modules instables :
{\small
\begin{equation}\label{T'}\tag{{\bf T'}}
0\rightarrow \omega_n L_n \rightarrow L_n \rightarrow \cdots \rightarrow \omega_{n-k+1}^{2^{k-1}-1}L_{n-k+1}\otimes 
\Sigma^{2^{k-1}-1}\F \xrightarrow{\delta_{k,n}} 
\omega_{n-k}^{2^{k}-1}L_{n-k}\otimes 
\Sigma^{2^{k}-1}\F\rightarrow \cdots \rightarrow \Sigma^{2^{n}-1}\F \rightarrow 0.
\end{equation}
}
D'apr\`es la proposition 4.2.1 de \cite{Tak99}, $\delta_{k,n}$ est donn\'e par 
$$
\delta_{k,n}\big(\omega_{n-k+1}^{2^{k-1}-1}\omega^{i_1,\cdots,i_{n-k+1}}\otimes \iota_{2^{k-1}-1}\big)
=\begin{cases} 0 & \text{si $i_{n-k+1}>1$,}\\
\omega_{n-k}^{2^{k}-1}\omega^{i_1,\cdots,i_{n-k}}\otimes \iota_{2^{k}-1} 
& \text{si $i_{n-k+1}=1$.}
\end{cases}
$$
On obtient donc un diagramme commutatif :
$$
\xymatrix{
\omega_{n-k+1}^{2^{k-1}-1}L_{n-k+1}\otimes \Sigma^{2^{k-1}-1}\F \ar[r]^-{\delta_{k,n}} \ar@{_(->}[d]
& \omega_{n-k}^{2^{k}-1}L_{n-k}\otimes \Sigma^{2^{k}-1}\F \ar@{_(->}[d]\\
L_{n-k+1}\otimes J(2^{k-1}-1)\ar[r]^-{f_{k,n}}	&	L_{n-k}\otimes J(2^{k}-1).
}$$
La suite \ref{T'}, dont la r\'ealisation g\'eom\'etrique est \ref{T}, 
d\'efinit un \'el\'ement dans $\Ext_{\Mc}^{n}(\Sigma^{2^n-1}\F,L'_n)$ qui \`a cause de la commutativit\'e ci-dessus 
est \'egal \`a $\nu_n$ de \ref{nu}. 
On \'etudiera ailleurs la r\'ealisation g\'eom\'etrique de la suite exacte de \ref{NST1}.

\subsection{La suite spectrale d'Adams}
On consid\`ere la suite spectrale d'Adams pour $\map(\L'(n),S^0)$ associ\'ee \`a $H^*(-,\Z/2)$  qui converge 
\`a $\pi_*(\map(\L'(n),\Sb^0))$ et dont le terme $E_2$ est 
$$E_2^{s,t}=\Ext_\A^{s,t}(\Z/2,L'_n)$$ et les diff\'erentielles $$d_r\colon E_r^{s,t}\rightarrow E_r^{s+r,t+r-1}.$$

\proofof{\ref{NST3}}
Pour $n\ge 2$, on doit v\'erifier que  
$$\pi^{k}(\L'(n))=
\begin{cases}\Z/2 &\text{si $k\in \{ n+2^n-1, n+2^{n-1}\}$,}\\
0 &\text{si $k\in [n+2^{n-2},+\infty)\setminus \{n+2^n-1, n+2^{n-1}\}$.}
\end{cases}
$$

D'apr\`es \ref{NST2}, la page $E_2$ de la suite spectrale pour 
$\pi^* (\L'(n))$ est repr\'esent\'ee comme dans la figure 1. 
On en d\'eduit le r\'esultat.

Pour $n=1$, on utilise la suite exacte longue 
d'homotopie de la cofibration $\Sb^1\rightarrow \L(1) \rightarrow \L'(1)$ et  
la conjecture de Segal \cite{AGM85} pour $\L(1)=\Sigma^\infty B\Z/2$, ce qui dit que 
$\pi^k(\L(1))=0$ si $k>0$ et $\pi^0(\L(1))=\Z_2$. On obtient alors
$\pi^1(\L'(1))=0$, 
$\pi^2(\L'(1))=\Z_2$ et $\pi^k(\L'(1))=0$ si $k>2$.
\qed

\begin{remark}
On peut composer l'application non triviale 
$\L'(n) \rightarrow \Sb^{n+2^n-1}$ avec l'inclusion de cellule de dimension minimale de $\L'(n)$, 
qui est $\Sb^{3\cdot 2^n-n-3}$. On obtient ainsi un \'el\'ement du groupe stable $\pi^s_{2^{n+1}-2n-2}$. 
Pour $n=1$, c'est l'application de degr\'e $2$, pour $n=2$ c'est $\eta^2$.
\end{remark}

\subsection{Sur la cohomotopie de $\M(n)_k$}
En fait le calcul de $\pi^*(\L'(n))$ en degré supérieur ou égal à $n+2^n-1$ 
peut être déduit directement des cofibrations 
de Takayasu et de la conjecture de Segal comme suit.

\begin{theorem}\label{NST4}
Pour tout $n\ge 1$ et tout $k\ge 1$, on a 
$$
\pi^t\M(n)_k=
\begin{cases} \pi^{2^{n-1}(k-1)+1}\M(1)_{2^{n-1}(k-1)+1} & 
\text{si $t=(2^n-1)(k-1)+n$,}\\
0 & \text{si $t>(2^n-1)(k-1)+n$.}
\end{cases}
$$
\end{theorem}

\begin{proof}
On pose $\alpha(n,k)=2^{n-1}(k-1)+1$ et $\beta(n,k)=(2^n-1)(k-1)+n$.
On fait une r\'ecurrence double sur les couples $(n,k)\ge (1,1)$ pour 
d\'emontrer la formule 
$$
\pi^t\M(n)_k=
\begin{cases} \pi^{\alpha(n,k)}\M(1)_{\alpha(n,k)} & 
\text{si $t=\beta(n,k)$,}\\
0 & \text{si $t>\beta(n,k)$.}
\end{cases}
$$

Supposons $k=1$. Il r\'esulte de la conjecture de Segal pour le spectre $\M(n)_1=\L(n)$ que 
$\pi^t\L(n)=0$ pour tout $t > 0$ et tout $n\ge 1$. 
La formule est v\'erifi\'ee.   

Supposons $n=1$. On fait une r\'ecurrence sur $k\ge 1$ pour montrer que
$\pi^t\M(1)_k=0$ si $t>\beta(1,k)=k$. Le cas $k=1$ est clair. Pour $k>1$, la cofibration 
$\Sb^{k-1} \rightarrow \M(1)_{k-1}\rightarrow \M(1)_k$ induit une suite exacte en cohomotopie : 
$$\pi^{t-1} \Sb^{k-1} \rightarrow \pi^t \M(1)_k \rightarrow \pi^t \M(1)_{k-1}.$$
Si $t> k$, on a $\pi^{t-1} \Sb^{k-1}=0$ et  $\pi^t \M(1)_{k-1}=0$ 
(d'apr\`es l'hypoth\`ese de r\'ecurrence), donc $\pi^t \M(1)_k=0$. On a ainsi v\'erifi\'e la formule 
pour $n=1$ et $k\ge 1$.

Supposons que $n\ge 2$, $k\ge 2$ et que la formule soit v\'erifi\'ee pour tout couple $(n',k')$ 
inf\'erieur au couple $(n,k)$ dans l'ordre lexicographique. 
La cofibration de Takayasu $\Sigma^{k-1} \M(n-1)_{2k-1} \rightarrow \M(n)_{k-1} \rightarrow \M(n)_k$ 
donne une suite exacte en cohomotopie :
$$\pi^{t-1}\M(n)_{k-1} \rightarrow \pi^{t-1} \Sigma^{k-1} \M(n-1)_{2k-1}
\rightarrow \pi^t \M(n)_k \rightarrow \pi^t \M(n)_{k-1}.$$
Si $t\ge \beta(n,k)$, par hypoth\`ese de r\'ecurrence double pour le couple $(n,k-1)$, on a 
$$\pi^{t-1}\M(n)_{k-1}=\pi^t \M(n)_{k-1}=0$$
car $$t>t-1 \ge \beta(n,k)-1 > \beta(n,k-1).$$
On en d\'eduit que, si $t\ge \beta(n,k)$, on a 
$$\pi^t \M(n)_k\cong \pi^{t-1} \Sigma^{k-1} \M(n-1)_{2k-1}\cong\pi^{t-k} \M(n-1)_{2k-1}.$$ 

Puisque $t\ge\beta(n,k) \Longleftrightarrow t-k \ge \beta(n-1,2k-1)$ et 
$\alpha(n,k)=\alpha(n-1,2k-1)$, 
l'hypoth\`ese de r\'ecurrence double pour le couple $(n-1,2k-1)$ donne 
la formule souhait\'ee pour $\pi^t\M(n)_k$.
\end{proof}

\begin{corollary} Soient $n\ge 1$ et $i\ge 0$. On a
$$
\pi^{\beta(n,2^i+1)} \M(n)_{2^i+1}=
\begin{cases} \Z_2 & \text{si $(n,i)=(1,0)$,}\\
\Z/2 & \text{sinon.}
\end{cases}
$$
\end{corollary}

\begin{proof}
On observe que $\alpha(n,2^i+1)=\alpha(n+i,2)$. D'apr\`es le th\'eor\`eme \ref{NST4}, on a alors
$$\pi^{\beta(n,2^i+1)}\M(n)_{2^i+1}=\pi^{\alpha(n,2^i+1)}\M(1)_{\alpha(n,2^i+1)}=
\pi^{\alpha(n+i,2)}\M(1)_{\alpha(n+i,2)}=\pi^{\beta(n+i,2)}\M(n+i)_2.$$
On applique le th\'eor\`eme \ref{NST3} au spectre $\L'(n+i)=\M(n+i)_2$ pour obtenir le r\'esultat.
\end{proof}

On a en particulier le r\'esultat suivant pour les espaces projectifs tronqu\'es  $\M(1)_{2^i+1}$:
\begin{corollary} \label{stunted}
On a
$$
\pi_S^{2^i+1} \M(1)_{2^i+1}=
\begin{cases} \Z_2 & \text{si $i=0$,}\\
\Z/2 & \text{si $i\ge 1$.}
\end{cases}
$$
\end{corollary}

\section{Appendice : Compl\'ements sur le facteur de Steinberg }\label{facteur-steinberg}

Pour $1\le k\le n$, 
l'image de l'idempotent $e_k\in \F[\gl_k]$ dans $\F[\gl_n]$ par 
l'inclusion canonique $\gl_k\hookrightarrow \gl_n$
utilisant les $m$ premi\`eres coordon\'ees de $\F^n$, par abus, 
se note aussi $e_k$.
Pour $1\le i \le n-1$, on d\'esigne par $e_{2,i}$ l'image de l'idempotent $e_2\in \F[\gl_2]$ dans 
$\F[\gl_n]$ par l'inclusion canonique $\gl_2\hookrightarrow \gl_n$ utilisant  
les  $i$-i\`eme et $(i+1)$-i\`eme coordon\'ees de $\F^n$. 

Dans \cite{Kuh87} 
Kuhn montre que la sous-alg\`ebre de $\F[\gl_n]$ engendr\'ee par $e_{2,1},\cdots,e_{2,n-1}$ 
est isomorphe \`a l'alg\`ebre de Hecke $\End_{\F[\gl_n]}(1_{B_n}^{\gl_n})$.  
En particulier, on a

\begin{proposition} [\cite{Kuh84,KP85}]\label{idempotent-relation}
\begin{enumerate}
\item $e_n$ est un produit de longueur maximale des 
$e_{2,1},\cdots,e_{2,n-1}$. 
\item $e_n=e_ne_{2,i}=e_{2,i}e_n$ pour tout $1\le i\le n-1$;
\item $e_n=e_{n-1}e_{2,n-1}e_{n-1}$.
\end{enumerate}
\end{proposition}

\proofof{\ref{cap-steinberg}} 
Il s'agit de montrer que 
$$L_n=e_n\cdot L_1^{\otimes n}=\omega_nM_n=\bigcap_{i=1}^{n-1} L_1^{\otimes i-1}\otimes L_2 \otimes L_1^{\otimes n-i-1}.$$

On commence par v\'erifier que $L_n=e_n\cdot L_1^{\otimes n}$ 
 en identifiant 
$M_1^{\otimes n}$ \`a l'alg\`ebre $\F[x_1,\cdots,x_n]$ et 
$L_1^{\otimes n}$ \`a l'id\'eal engendr\'e par $x_1\cdots x_n$. 
On observe d'abord que $e_n\cdot L_1^{\otimes n}$ est en fait facteur direct de 
$L_1^{\otimes n}$ puisque l'idempotent $e_n$ induit un endomorphisme de $L_1^{\otimes n}$. 
A cet effet, on se sert de \ref{idempotent-relation} (1), ce qui dit que 
$e_n$ est certain produit des $e_{2,i}$, pour se ramener au cas 
$n=2$, ce qui est ais\'ement v\'erifi\'e par calcul direct.
Comme $ L_1^{\otimes n}$ est facteur direct de $M_1^{\otimes n}$, $e_n\cdot L_1^{\otimes n}$ est  
facteur direct de $M_n=e_n\cdot M_1^{\otimes n}$. La rigidit\'e de l'isomorphisme $M_n\cong L_n\oplus L_{n-1}$ 
implique alors que $L_n=e_n\cdot L_1^{\otimes n}$. En effet, 
si l'on note $I_{n-1}\in M_n(\F)$ la matrice $\diag(1,\cdots,1,0)$,
Kuhn  a montr\'e dans \cite{Kuh87} que $\End_\U(M_n)$ est l'espace vectoriel engendr\'e par les idempotents 
$e_n$ et $e_nI_{n-1}e_n$. 
En particulier, le facteur direct $L_{n-1}$ de $M_n$ correspond \`a l'idempotent $e_nI_{n-1}e_n$. 
Comme l'action de $e_nI_{n-1}$ est triviale sur $e_n\cdot L_1^{\otimes n}\subset L_1^{\otimes n}$, 
on  obtient  $L_n=e_n\cdot L_1^{\otimes n}$.

On v\'erifie ensuite que 
$\omega_ne_n\cdot M_1^{\otimes n}=e_n\cdot L_1^{\otimes n}$. 
Comme l'invariant $\omega_n$ est divisible par $x_1\cdots x_n$, il est clair que 
$\omega_ne_n\cdot M_1^{\otimes n}\subset e_n\cdot L_1^{\otimes n}$. 
Soit $f\in e_n\cdot L_1^{\otimes n}\subset L_1^{\otimes n}$. Alors $f$ est $B_n$-invariant et
divisible par $x_1\cdots x_n$. Une observation \'el\'ementaire, dûe \`a H. M\`ui, dit que 
si $f$ est $B_n$-invariant et divisible par $x_k$, alors $f$ est \'egalement divisible 
par 
$$V_k=\prod_{\lambda_i\in \F}(\lambda_1x_1+\cdots+\lambda_{k-1}x_{k-1}+x_k).$$ On en d\'eduit que 
$f$ est divisible par $\omega_n=V_1\cdots V_n$. Soit $f=\omega_nf'$. On a donc 
$f=e_n\cdot f=\omega_n e_n\cdot f'$ est un \'el\'ement de $\omega_ne_n\cdot M_1^{\otimes n}$ et 
ainsi obtient l'inclusion $e_n\cdot L_1^{\otimes n}\subset \omega_ne_n\cdot M_1^{\otimes n}$. 
D'o\`u $e_n\cdot L_1^{\otimes n}=\omega_nM_n$.

Enfin l'identification 
$$e_n\cdot L_1^{\otimes n}=\bigcap_{i=1}^{n-1} L_1^{\otimes i-1}\otimes L_2 \otimes L_1^{\otimes n-i-1},$$ 
dûe \`a Kuhn \cite{Kuh84}, est d\'emontr\'ee en utilisant \ref{idempotent-relation} (1) et (2). 
Si $x$ est un \'el\'ement de $e_n\cdot L_1^{\otimes n}$, alors $x$ est invariant par $e_n$. 
Comme $e_{2,i}e_n=e_n$ pour tout $1\le i\le n-1$, 
$x$ est aussi invariant par $e_{2,i}$ pour tout $1\le i\le n-1$. D'o\`u $x$ appartient \`a
$\bigcap_{i=1}^{n-1} L_1^{\otimes i-1}\otimes L_2 \otimes L_1^{\otimes n-i-1}$. Si $x$ est un \'el\'ement de 
$\bigcap_{i=1}^{n-1} L_1^{\otimes i-1}\otimes L_2 \otimes L_1^{\otimes n-i-1}$, alors $x$ est invariant par $e_n$ 
car $e_n$ est certain produit des $e_{2,i}$, d'o\`u $x\in e_n\cdot L_1^{\otimes n}$.
\qed

L'alg\`ebre de Dickson $D(k)$
est la sous-alg\`ebre des invariants sous l'action du groupe
lin\'eaire $\gl_k$ sur $\F[x_1,\dots,x_k]$. Soit $\omega_k=\det (x_{j}^{2^{i-1}})_{1\le i,j\le k}$
la classe de Dickson sup\'erieure de $D(k)$.

\begin{proposition}\label{diagonal} Si $i_1> 2i_2>\cdots> 2^{n-1}i_{n}>0$, on a
$$e_{n}\cdot \omega_1^{i_1-2i_2}\cdots
\omega_{n-1}^{i_{n-1}-2i_{n}} \omega_{n}^{i_{n}}= (e_{n-1}\cdot \omega_1^{i_1-2i_2}\cdots
\omega_{n-2}^{i_{n-2}-2i_{n-1}} \omega_{n-1}^{i_{n-1}})\cdot x_{n}^{i_{n}}+
\sum_i f_i\cdot x_{n}^i,$$
pour certains $i>i_{n}$ et $f_i\in L_{n-1}\subset \F[x_1,\ldots,x_{n-1}]$. 
\end{proposition}

\begin{proof}
Notons d'abord que $\omega_{n}$ est $\gl_n$-invariant et l'on a un d\'eveloppement 
$$\omega_{n}^{i_n}=\omega_{n-1}^{2i_n}\cdot x_n^{i_n}+\sum_{i>i_n}g_i\cdot x_n^i$$
pour certains $g_i\in \F[x_1,\cdots,x_{n-1}]$. Il suffit donc de montrer que 
$$e_n\cdot \omega_{1}^{i_1}\cdots \omega_{n-1}^{i_{n-1}}=
e_{n-1}\cdot \omega_{1}^{i_1}\cdots \omega_{n-1}^{i_{n-1}}+\sum_{j>0}h_j\cdot x_n^j$$
pour certains $h_j\in \F[x_1,\cdots,x_{n-1}]$. De mani\`ere \'equivalente, il suffit de montrer 
que 
$$I_{n-1}e_n\cdot \omega_{1}^{i_1}\cdots \omega_{n-1}^{i_{n-1}}=
e_{n-1}\cdot \omega_{1,0}^{i_1}\cdots \omega_{n-1}^{i_{n-1}}$$
avec $I_{n-1}=\diag(1,\cdots,1,0)\in M_n(\F)$.
Posons $Q:=Q(x_1,\cdots,x_{n-1})=\omega_{1}^{i_1}\cdots \omega_{n-1}^{i_{n-1}}$. On a 
\begin{eqnarray*}
I_{n-1}e_n\cdot Q(x_1,\cdots,x_{n-1})&=&I_{n-1}e_nI_{n-1}\cdot Q(x_1,\cdots,x_{n-1})\\
				&=& I_{n-1}\big (e_{n-1}e_{2,n-1}e_{n-1}\big )I_{n-1}\cdot Q(x_1,\cdots,x_{n-1})\\
				&=& e_{n-1}\big (I_{n-1}e_{2,n-1}I_{n-1}\big )e_{n-1}\cdot Q(x_1,\cdots,x_{n-1}).
\end{eqnarray*}
Comme  $I_1e_2I_1=\diag(1,0)+\diag(0,0)$, on obtient  
\begin{align*}I_{n-1}e_n\cdot Q
&=e_{n-1}\diag(1,\cdots,1,1,0)e_{n-1}\cdot Q+e_{n-1}\diag(1,\cdots,1,0,0)e_{n-1}\cdot Q\\
&=e_{n-1}\diag(1,\cdots,1,1,0)e_{n-1}\cdot Q	
\quad \quad \text{(comme $e_{n-1}\cdot Q$ est divisible par $x_{n-1}$)}\\
&=e_{n-1}e_{n-1}\cdot Q\\
&=e_{n-1}\cdot Q.
\end{align*}
La proposition suit.
\end{proof}

\proofof{\ref{base-steinberg}} 
Il suffit de v\'erifier que
\begin{equation} \label{dickson}
\Sq^{i_1+1}\cdots \Sq^{i_n+1}(\frac{1}{x_1\cdots x_n})
={\bar{\Sigma}_n}\cdot \omega_1^{i_1-2i_2}\cdots \omega_{n-1}^{i_{n-1}-2i_n}
\omega_n^{i_n}
\end{equation}
pour $i_1> 2i_2>\cdots >2^{n-1}i_n\geq 0$.
On fait une r\'ecurrence sur $n$. Le cas $n=1$ est trivial :
$$\Sq^{i_1+1}(\frac{1}{x_1})=x_1^{i_1}.$$
Supposons que la proposition soit vraie pour tous les entiers inf\'erieurs \`a $n$.
Posons $J:=(i_2+1,\cdots,i_n+1)$. On a
\begin{equation*}
\begin{split}
\Sq^J&(\frac{1}{x_1\cdots x_n})=P_J(x_1,\ldots,x_n)+
\sum_{1\le i\le n}\frac{1}{x_i}\Sq^J(\frac{1}{x_1\cdots \hat{x}_i\cdots x_n})\\
&+\sum_{1\le i<j\le n}\frac{1}{x_ix_j}
\Sq^J(\frac{1}{x_1\cdots \hat{x}_i\cdots \hat{x}_j\cdots x_n})+
\sum_{1\le i<j<k\le n}\frac{1}{x_ix_jx_k}
\Sq^J(\frac{1}{x_1\cdots \hat{x}_i\cdots \hat{x}_j\cdots \hat{x}_k\cdots x_n})+\cdots,
\end{split}
\end{equation*}
ou $P_J(x_1,\cdots,x_n)\in \F[x_1,\ldots,x_n]$ est un polyn\^ome de degr\'e $i_2+\cdots+i_n-1$.
Par hypoth\`ese de r\'ecurrence et instabilit\'e,
on voit que les termes de la deuxi\`eme ligne sous-d\'essus sont nuls. De même,
on a
$$\Sq^{i_1+1}P_J(x_1,\cdots,x_n)=0$$ par instabilit\'e. On obtient donc
$$\Sq^I(\frac{1}{x_1\cdots x_n})=
\Sq^{i_1+1}(\sum_{1\le i\le n}\frac{1}{x_i}\Sq^J(\frac{1}{x_1\cdots \hat{x}_i\cdots x_n})).$$
L'\'egalit\'e (\ref{dickson}) r\'esulte de l'hypoth\`ese de r\'ecurrence et du lemme suivant.
\begin{lemma}
\begin{equation*}
\begin{split}
\Sq^{i_1+1}&(\frac{\omega_1^{i_2-2i_3}(x_2)\cdots \omega_{n-2}^{i_{n-1}-2i_n}(x_2,\ldots,x_{n-1})
\omega_{n-1}^{i_n}(x_2,\ldots,x_n)}{x_1})\\
&=\omega_1^{i_1-2i_2}(x_1)\cdots
\omega_{n-1}^{i_{n-1}-2i_n}(x_1,\ldots,x_{n})\omega_{n}^{i_n}(x_1,\ldots,x_n).
\end{split}
\end{equation*}
\end{lemma}
On va montrer le lemme en utilisant le {\it carr\'e total stable} d\'efini comme suit.
Soit $M$ un $\A$-module. Supposons que $H^*(B\Z/2)\cong \F[x]$ avec $|x|=1$. Le
carr\'e total stable $\Sc\colon M\rightarrow \F[x,x^{-1}]\hat{\otimes}M$
est donn\'e par $$\Sc(z)=\sum_{i\geq 0}x^{-i}\otimes \Sq^i(z), \quad z\in M.$$
Notons que $\St(z)=x^{|z|}\Sc(z)$ est le {\it carr\'e total instable} de $z$.
De plus $\Sc$ est multiplicatif si $M$ est une $\A$-alg\`ebre. 
H. M\`ui \cite{Mui75} a montr\'e que $$\St(V_n(x_1,\ldots,x_n))=V_{n+1}(x,x_1,\ldots,x_n).$$
Ici
$$V_i(x_1,\ldots,x_i):=\prod_{(a_1,\ldots,a_{i-1})\in\F^{i-1}}(a_1x_1+\cdots+a_{i-1}x_{i-1}+x_i)$$
et
$$\F[x_1,\ldots,x_n]^{B_n}\cong \F[V_1,\ldots,V_n].$$
Utilisant la relation $\omega_k=V_1\cdots V_k$, on a
$$\Sc\bigg(\omega_k(x_2,\cdots,x_{k+1})\bigg)=\frac{\omega_{k+1}(x,x_2,\cdots,x_{k+1})}{x^{2^k}}.$$
Posons
$$z=\frac{\omega_1^{i_2-2i_3}(x_2)\cdots \omega_{n-2}^{i_{n-1}-2i_n}(x_2,\ldots,x_{n-1})\omega_{n-1}^{i_n}(x_2,\ldots,x_n)}{x_1}.$$
On obtient
\begin{equation*}
\Sc (z)=\frac{\omega_1^{i_1-2i_2}(x)\cdots \omega_{n-1}^{i_{n-1}-2i_n}(x,x_2,\ldots,x_{n})
\omega_{n}^{i_n}(x,x_2,\ldots,x_n)}{x^{i_1}}\Sc(\frac{1}{x_1})=:\frac{\omega_I(x)}{x^{i_1}}\Sc (\frac{1}{x_1}).
\end{equation*}
L'admissibilit\'e de la suite $I$ permet d'\'ecrire $\omega_I(x)$ sous la forme 
$$\omega_I(x)=\sum_{j=0}^mf_jx^j$$
avec $f_j\in \F[x_2,\cdots,x_n]$. D'autre part, l'action de l'alg\`ebre de Steenrod sur $\frac{1}{x_1}$ donne 
$$\Sc (\frac{1}{x_1})=\sum_{i\ge 0}\frac{x_1^{i-1}}{x^i}.$$
Il suit 
$$\Sc(z)=\frac{1}{x^{i_1+1}}\sum_{i\ge 0}\sum_{0\le j\le m}f_jx_1^ix^{j-i}.$$
On voit clairement que le coefficient de $\frac{1}{x^{i_1+1}}$ dans la s\'erie formelle
$\Sc(z)$ est $\sum_{j=0}^mf_jx_1^j=\omega_I(x_1)$. Le lemme est d\'emontr\'e.
\qed

\bibliographystyle{amsalpha}

\providecommand{\bysame}{\leavevmode\hbox to3em{\hrulefill}\thinspace}
\providecommand{\MR}{\relax\ifhmode\unskip\space\fi MR }
\providecommand{\MRhref}[2]{%
  \href{http://www.ams.org/mathscinet-getitem?mr=#1}{#2}
}
\providecommand{\href}[2]{#2}

\bigskip
NGUYEN DANG HO HAI

Universit\'e de Hu\'e, Coll\`ege des Sciences, 77 Rue Nguyen Hue, Hue Ville, VIETNAM

et

LAGA, UMR 7539 du CNRS, Universit\'e Paris 13

99, Av. J-B Cl\'ement, 93430 Villetaneuse, FRANCE

nguyen@math.univ-paris13.fr

\bigskip

LIONEL SCHWARTZ

LAGA, UMR 7539 du CNRS, Universit\'e Paris 13

99, Av. J-B Cl\'ement, 93430 Villetaneuse, FRANCE

schwartz@math.univ-paris13.fr

\bigskip

TRAN NGOC NAM

Universit\'e Nationale du Vietnam, Coll\`ege des Sciences,  334 Rue Nguyen Trai, Hanoi, VIETNAM

bruce-nam@hotmail.com

\end{document}